# ANALYSIS OF SPDES ARISING IN PATH SAMPLING PART II: THE NONLINEAR CASE

BY M. HAIRER[1], A. M. STUART[1,2] AND J. VOSS[2]

*University of Warwick*

In many applications, it is important to be able to sample paths of SDEs conditional on observations of various kinds. This paper studies SPDEs which solve such sampling problems. The SPDE may be viewed as an infinite-dimensional analogue of the Langevin equation used in finite-dimensional sampling. In this paper, conditioned nonlinear SDEs, leading to nonlinear SPDEs for the sampling, are studied. In addition, a class of preconditioned SPDEs is studied, found by applying a Green's operator to the SPDE in such a way that the invariant measure remains unchanged; such infinite dimensional evolution equations are important for the development of practical algorithms for sampling infinite dimensional problems.

The resulting SPDEs provide several significant challenges in the theory of SPDEs. The two primary ones are the presence of nonlinear boundary conditions, involving first order derivatives, and a loss of the smoothing property in the case of the pre-conditioned SPDEs. These challenges are overcome and a theory of existence, uniqueness and ergodicity is developed in sufficient generality to subsume the sampling problems of interest to us. The Gaussian theory developed in Part I of this paper considers Gaussian SDEs, leading to linear Gaussian SPDEs for sampling. This Gaussian theory is used as the basis for deriving nonlinear SPDEs which affect the desired sampling in the nonlinear case, via a change of measure.

**1. Introduction.** The purpose of this paper is to provide rigorous justification for a recently introduced stochastic partial differential equation (SPDE)-based approach to infinite dimensional sampling problems [14, 22]. The methodology has been developed to solve a number of sampling problems arising from stochastic differential equations (SDEs—assumed to be finite-dimensional unless stated otherwise), conditional on observations.

Received May 2006; revised March 2007.
[1]Supported by EPSRC Grant EP/E002269/1.
[2]Supported by EPSRC and ONR Grant N00014-05-1-0791.

The setup is as follows. Consider the SDE

$$(1.1) \qquad dX = AX\, du + f(X)\, du + B\, dW^x, \qquad X(0) = x^-,$$

where $f(x) = -BB^*\nabla V(x)$, $V: \mathbb{R}^d \to \mathbb{R}$, $B \in \mathbb{R}^{d \times d}$ is invertible and $W^x$ is a standard $d$-dimensional Brownian motion. We consider three sampling problems associated with (1.1):

1. *free path sampling*, to sample paths of (1.1) unconditionally;

2. *bridge path sampling*, to sample paths of (1.1) *conditional* on knowing $X(1) = x^+$;

3. *nonlinear filter/smoother*, to sample paths of (1.1), conditional on knowledge of $(Y(u))_{u \in [0,1]}$ solving

$$(1.2) \qquad dY = \tilde{A}X\, dt + \tilde{B}\, dW^y, \qquad Y(0) = 0,$$

where $\tilde{A} \in \mathbb{R}^{m \times d}$ is arbitrary, $\tilde{B} \in \mathbb{R}^{m \times m}$ is invertible and $W^y$ is a standard $m$-dimensional Brownian motion.

The methodology proposed in [22] is to extend the finite dimensional Langevin sampling technique [21] to infinite-dimensional problems such as those listed above as 1 to 3. This leads to SPDEs which are ergodic and have stationary measure which solves the desired sampling problem.

We believe that an infinite dimensional sampling technique can be derived by taking the (formal) density of the target distribution and mimicking the procedure from the finite dimension Langevin method. In this paper, we provide a rigorous justification for this claim in the case of equation (1.1), where the drift is linear plus a gradient, the noise is additive and, in case 3, observations arise linearly, as in (1.2). A conjecture for the case of general drift, and for a nonlinear observation equation in place of (1.2), is described in Section 9 at the end of the paper.

For the problems considered here, the resulting SDPEs are of the form

$$(1.3) \qquad dx = (BB^*)^{-1}\partial_u^2 x\, dt - \nabla\Phi(x)\, dt + \sqrt{2}\, dw(t),$$

and generalizations, where $w$ is a cylindrical Wiener process (so that $\frac{\partial w}{\partial t}$ is space-time white noise) and $\Phi$ is some real-valued function on $\mathbb{R}^d$. [Note that the "potential" $\Phi$ is different from $V$; see (5.3) below.] For problem 1, the resulting SPDE is not a useful algorithmic framework in practice as it is straightforward to generate unconditioned, independent samples from 1 by application of numerical methods for SDE initial value problems [15]; the Langevin method generates correlated samples and, hence, has larger variance in any resulting estimators. However, we include analysis of problem 1 because it contributes to the understanding of subsequent SPDE-based approaches. For problems 2 and 3, we believe that the proposed methodology is, potentially, the basis for efficient Markov chain Monte Carlo (MCMC)-based sampling techniques. Some results about how such an MCMC method could be implemented in practice can be found in [1] and [13].

The resulting MCMC method, when applied to problem 3, results in a new method for solving nonlinear filtering/smoothing problems. This method differs substantially from traditional methods like particle filters which are based on the Zakai equation. While the latter equation describes the density of the conditional distribution of the signal at fixed times $t$, our proposed method samples full paths from the conditional distribution; statistical quantities can then be obtained by considering ergodic averages. Consequently, while the proposed method cannot easily be applied in online situations, it provides dynamic information on the paths, which cannot be so easily read off the solutions of the Zakai equation. Another difference is that the independent variables in the Zakai equation are in $\mathbb{R}^d$, whereas equation (1.3) is always indexed by $[0,\infty) \times [0,1]$ and only takes values in $\mathbb{R}^d$. Thus, the proposed method should be advantageous in high dimensions. For further discussion and applications, see [13].

In making such methods as efficient as possible, we are lifting ideas from finite-dimensional Langevin sampling into our infinite-dimensional situation. One such method is to use preconditioning, which changes the evolution equation, whilst preserving the stationary measure, in an attempt to roughly equalize relaxation rates in all modes. This leads to SPDEs of the form

$$(1.4) \qquad dx = \mathcal{G}(BB^*)^{-1}\partial_u^2 x \, dt - \mathcal{G}\nabla\Phi(x) \, dt + \sqrt{2}\mathcal{G}^{1/2} \, dw(t),$$

and generalizations, where $w$, again, is a cylindrical Wiener process. In the finite-dimensional case, it is well known that the invariant measure for (1.4) is the same as for (1.3). In this paper, we will study the methods proposed in [1] which precondition the resulting infinite-dimensional evolution equation (1.3) by choosing $\mathcal{G}$ as a Green's operator. We show that equation (1.4), in its stationary measure, still samples from the desired distribution.

For both preconditioned and unpreconditioned equations, the analysis leads to mathematical challenges. First, when we are not conditioning on the endpoint (in problems 1 and 3), we obtain an SPDE with a nonlinear boundary condition of the form

$$\partial_u x(t,1) = f(x(t,1)) \qquad \forall t \in (0,\infty),$$

where $f$ is the drift of the SDE (1.1). In the abstract formulation using Hilbert-space-valued equations, this translates into an additional drift term of the form $f(x(t,1))\delta_1$, where $\delta_1$ is the delta distribution at $u=1$. This forces us to consider equations with values in the Banach space of continuous functions (so that we can evaluate the solution $x$ at the point $u=1$) and to allow the drift to take distributions as values. Unfortunately, the theory for this situation is not well developed in the literature. Therefore, we here provide proofs for the existence and uniqueness of solutions for the SPDEs considered. This machinery is not required for problem 2, as the Hilbert space setting [4, 5, 6, 25] can be used there.

We also prove ergodicity of these SPDEs. Here, a second challenge comes from the fact that we consider the preconditioned equation (1.4). Since we want to precondition with operators $\mathcal{G}$ which are close to $(\partial_u^2)^{-1}$, it is not possible to use smoothing properties of the heat semigroup anymore and the resulting process no longer has the strong Feller property. Instead, we show that the process has the recently introduced asymptotic strong Feller property (see [12]) and use this to show existence of a unique stationary measure for the preconditioned case.

The paper is split into two parts. The first part, consisting of Sections 2, 3 and 4, introduces the general framework, while the second part, starting at Section 5, uses this framework to solve the three sampling problems stated above. Readers only interested in the applications can safely skip the first part on first reading. The topics presented there are mainly required to understand the proofs in the second part.

The two parts are organized as follows. In Section 2, we introduce the technical framework required to give sense to equations like (1.3) and (1.4) as Hilbert-space-valued SDEs. The main results of this section are Theorem 3.4 and 3.6, showing the global existence of solutions of these SDEs. In Section 3, we identify a stationary distribution of these equations. This result is a generalization of a result by Zabczyk [25]; the generalization allows us to consider the Banach-space-valued setting required for the nonlinear boundary conditions and is also extended to consider the preconditioned SPDEs. In Section 4, we show that the stationary distribution is unique and that the considered equations are ergodic (see Theorems 4.10 and 4.11). This justifies their use as the basis for an MCMC method.

In the second part of the paper, we apply the abstract theory to derive SPDEs which sample conditioned diffusions. Section 5 outlines the methodology. Then, in Sections 6, 7 and 8, we discuss the sampling problems 1, 2 and 3, respectively, proving the desired property for both the SPDEs proposed in [22] and the preconditioned method proposed in [1]. In the case 2, bridges, the SPDE whose invariant measure is the bridge measure was also derived in one dimension in [20]. In Section 9, we give a heuristic method to derive SPDEs for sampling, which applies in greater generality than the specific setups considered here. Specifically, we show how to derive the SPDE when the drift vector field in (1.1) is not of the form "linear plus gradient"; for signal processing, we show how to extend beyond the case of linear observation equation (1.2). This section will be of particular interest to the reader concerned with *applying* the technique for sampling which we study here. The gap between what we conjecture to be the correct SPDEs for sampling in general and the cases considered here points to a variety of open and interesting questions in stochastic analysis; we highlight these.

To avoid confusion, we use the following naming convention. Solutions to SDEs like (1.1), which give our target distributions, are denoted by upper case letters like $X$. Solutions to infinite-dimensional Langevin equations like (1.3), which we use to sample from these target distributions, are denoted by lower case letters like $x$. The variable which is time in equation (1.1) and space in (1.3) is denoted by $u$ and the time direction of our infinite-dimensional equations, which indexes our samples, is denoted by $t$.

**2. The abstract framework.** In this section, we introduce the abstract setting for our Langevin equations, proving existence and uniqueness of global solutions. We treat the nonpreconditioned equation (1.3) and the preconditioned equation (1.4) separately. The two main results are Theorems 2.6 and 2.10. Both cases will be described by stochastic evolution equations taking values in a real Banach space $E$ continuously embedded into a real separable Hilbert space $\mathcal{H}$. In our applications in the later sections, the space $\mathcal{H}$ will always be the space of $L^2$ functions from $[0,1]$ to $\mathbb{R}^d$ and $E$ will be some subspace of the space of continuous functions.

Our application requires the drift to be a map from $E$ to $E^*$. This is different from the standard setup as found in, for example, [5], where the drift is assumed to take values in the Hilbert space $\mathcal{H}$.

2.1. *The nonpreconditioned case.* In this subsection, we consider semilinear SPDEs of the form

$$(2.1) \qquad dx = \mathcal{L}x\,dt + F(x)\,dt + \sqrt{2}\,dw(t), \qquad x(0) = x_0,$$

where $\mathcal{L}$ is a linear operator on $\mathcal{H}$, the drift $F$ maps $E$ into $E^*$, $w$ is a cylindrical Wiener process on $\mathcal{H}$ and the process $x$ takes values in $E$. We seek a mild solution of this equation, defined precisely below.

Recall that a closed, densely defined operator $\mathcal{L}$ on a Hilbert space $\mathcal{H}$ is called *strictly dissipative* if there exists $c > 0$ such that $\langle x, \mathcal{L}x \rangle \leq -c\|x\|^2$ for every $x \in \mathcal{D}(\mathcal{L})$. We make the following assumptions on $\mathcal{L}$.

(A1) Let $\mathcal{L}$ be a self-adjoint, strictly dissipative operator on $\mathcal{H}$ which generates an analytic semigroup $S(t)$. Assume that $S(t)$ can be restricted to a $C_0$-semigroup of contraction operators on $E$.

Since $-\mathcal{L}$ is self-adjoint and positive, one can define arbitrary powers of $-\mathcal{L}$. For $\alpha \geq 0$, let $\mathcal{H}^\alpha$ denote the domain of the operator $(-\mathcal{L})^\alpha$ endowed with the inner product $\langle x, y \rangle_\alpha = \langle (-\mathcal{L})^\alpha x, (-\mathcal{L})^\alpha y \rangle$. We further define $\mathcal{H}^{-\alpha}$ as the dual of $\mathcal{H}^\alpha$ with respect to the inner production $\mathcal{H}$ (so that $\mathcal{H}$ can be seen as a subspace of $\mathcal{H}^{-\alpha}$). Denote the Gaussian measure with mean $\mu \in \mathcal{H}$ and covariance operator $\mathcal{C}$ on $\mathcal{H}$ by $\mathcal{N}(\mu, \mathcal{C})$.

(A2) There exists an $\alpha \in (0, 1/2)$ such that $\mathcal{H}^\alpha \subset E$ (densely), $(-\mathcal{L})^{-2\alpha}$ is nuclear in $\mathcal{H}$ and the Gaussian measure $\mathcal{N}(0, (-\mathcal{L})^{-2\alpha})$ is concentrated on $E$.

This condition implies that the stationary distribution $\mathcal{N}(0, (-\mathcal{L})^{-1})$ of the linear equation

$$dz = \mathcal{L}z\, dt + \sqrt{2}\, dw(t) \tag{2.2}$$

is concentrated on $E$.

Under assumption (A2), we have the following chain of inclusions:

$$\mathcal{H}^{1/2} \hookrightarrow \mathcal{H}^\alpha \hookrightarrow E \hookrightarrow \mathcal{H} \hookrightarrow E^* \hookrightarrow \mathcal{H}^{-\alpha} \hookrightarrow \mathcal{H}^{-1/2}.$$

Since we assumed that $E$ is continuously embedded into $\mathcal{H}$, each of the corresponding inclusion maps is bounded and continuous. Therefore, we can, for example, find a constant $c$ with $\|x\|_{E^*} \leq c\|x\|_E$ for all $x \in E$. Later, we will use the fact that, in this situation, there exist constants $c_1$ and $c_2$ with

$$\|S(t)\|_{E^* \to E} \leq c_1 \|S(t)\|_{\mathcal{H}^{-\alpha} \to \mathcal{H}^\alpha} \leq c_2 t^{-2\alpha}. \tag{2.3}$$

We begin our study of equation (2.1) with the following, preliminary result which shows that the linear equation takes values in $E$.

LEMMA 2.1. *Assume* (A1) *and* (A2) *and define the $\mathcal{H}$-valued process $z$ by the stochastic convolution*

$$z(t) = \sqrt{2} \int_0^t S(t-s)\, dw(s) \qquad \forall t \geq 0, \tag{2.4}$$

*where $w$ is a cylindrical Wiener process on $\mathcal{H}$. Then, $z$ has an $E$-valued continuous version. Furthermore, its sample paths are almost surely $\beta$-Hölder continuous for every $\beta < 1/2 - \alpha$. In particular, for such $\beta$, there exist constants $C_{p,\beta}$ such that*

$$\mathbb{E} \sup_{s \leq t} \|z(s)\|_E^p \leq C_{p,\beta} t^{\beta p} \tag{2.5}$$

*for every $t \leq 1$ and every $p \geq 1$.*

PROOF. Let $i$ be the inclusion map from $\mathcal{H}^\alpha$ into $E$ and $j$ be the inclusion map from $E$ into $\mathcal{H}$. Since $\langle x, y \rangle_\alpha = \langle (-\mathcal{L})^{2\alpha} jix, jiy \rangle$ for every $x, y$ in $\mathcal{H}^\alpha$, one has

$$\langle x, jiy \rangle = \langle i^* j^* x, y \rangle_\alpha = \langle (-\mathcal{L})^{2\alpha} jii^* j^* x, jiy \rangle$$

for every $x \in \mathcal{H}$ and every $y \in \mathcal{H}^\alpha$. Since $\mathcal{H}^\alpha$ is dense in $\mathcal{H}$, this implies that $jii^*j^* = (-\mathcal{L})^{-2\alpha}$. Thus, (A2) implies that $ii^*$ is the covariance of a Gaussian measure on $E$, which is sometimes expressed by saying that the map $i$ is $\gamma$-*radonifying*.

The first part of the result then follows directly from [3], Theorem 6.1. Conditions (i) and (ii) there are direct consequences of our assumptions (A1) and (A2). Condition (iii) there states that the reproducing kernel Hilbert space $\mathcal{H}_t$ associated with the bilinear form $\langle x, \mathcal{L}^{-1}(e^{\mathcal{L}t} - 1)y\rangle$ has the property that the inclusion map $\mathcal{H}_t \to E$ is $\gamma$-radonifying. Since we assumed that $\mathcal{L}$ is strictly dissipative, it follows that $\mathcal{H}_t = \mathcal{H}^{1/2}$. Since we just checked that the inclusion map from $\mathcal{H}^{1/2}$ into $E$ is $\gamma$-radonifying, the required conditions hold.

If we can show that $\mathbb{E}\|z(t+h) - z(t)\|_E \leq C|h|^{1/2-\alpha}$ for some constant $C$ and for $h \in [0,1]$, then the second part of the result follows from Fernique's theorem [10] combined with Kolmogorov's continuity criterion [19], Theorem 2.1. One has

$$\mathbb{E}\|z(t+h) - z(t)\|_E \leq \mathbb{E}\|S(h)z(t) - z(t)\|_E + \sqrt{2}\mathbb{E}\left\|\int_0^h S(s)\,dw(s)\right\|_E = T_1 + T_2.$$

The random variable $z(t)$ is Gaussian on $\mathcal{H}$ with covariance given by

$$Q_t = (-\mathcal{L})^{-1}(I - S(2t)).$$

This shows that the covariance of $(S(h) - I)z(t)$ is given by

$$(S(h) - I)Q_t(S(h) - I) = (-\mathcal{L})^{-\alpha}A_\alpha(I - S(2t))A_\alpha(-\mathcal{L})^{-\alpha},$$

with

$$A_\alpha = (-\mathcal{L})^{\alpha-1/2}(S(h) - I).$$

Since (A2) implies that $(-\mathcal{L})^{-\alpha}$ is $\gamma$-radonifying from $\mathcal{H}$ to $E$ and $(S(2t) - I)$ is bounded by 2 as an operator from $\mathcal{H}$ to $\mathcal{H}$, we have

$$T_1 \leq C\|A_\alpha\|_{L(\mathcal{H})} \leq C|h|^{1/2-\alpha},$$

where the last inequality follows from the fact that $\mathcal{L}$ is self-adjoint and strictly dissipative. The bound on $T_2$ can be obtained in a similar way. From Kolmogorov's continuity criterion, we get that $z$ has a modification which is $\beta$-Hölder continuous for every $\beta < 1/2 - \alpha$.

Since we now know that $z$ is Hölder continuous, the expression

(2.6) $$\sup_{\substack{s,t\in[0,1]\\s\neq t}} \frac{\|z(s) - z(t)\|_E}{|t-s|^\beta}$$

is finite almost surely. Since the field $\frac{z(s)-z(t)}{|t-s|^\beta}$ is Gaussian, it then follows from Fernique's theorem that (2.6) also has moments of every order. $\square$

REMARK 2.2. The standard factorization technique ([5], Theorem 5.9) does not apply in this situation since, in general, there exists no interpolation space $\mathcal{H}^\beta$ such that $\mathcal{H}^\beta \subset E$ and $z$ takes values in $\mathcal{H}^\beta$: for $\mathcal{H}^\beta \subseteq E$, one would require $\beta > 1/4$, but the process takes values in $\mathcal{H}^\beta$ only for $\beta < 1/4$. Lemma 2.1 should rather be considered as a slight generalization of [5], Theorem 5.20.

DEFINITION 2.3. The *subdifferential* of the norm $\|\cdot\|_E$ at $x \in E$ is defined as

$$\partial \|x\|_E = \{x^* \in E^* | x^*(x) = \|x\|_E \text{ and } x^*(y) \leq \|y\|_E \; \forall y \in E\}.$$

This definition is equivalent to the one in [5], Appendix D and, by the Hahn–Banach theorem, the set $\partial \|x\|_E$ is nonempty. We use the subdifferential of the norm to formulate the conditions on the nonlinearity $F$. Here and below, $C$ and $N$ denote arbitrary positive constants that may change from one equation to the next.

(A3) The nonlinearity $F: E \to E^*$ is Fréchet differentiable with

$$\|F(x)\|_{E^*} \leq C(1 + \|x\|_E)^N \quad \text{and} \quad \|DF(x)\|_{E \to E^*} \leq C(1 + \|x\|_E)^N$$

for every $x \in E$.

(A4) There exists a sequence of Fréchet differentiable functions $F_n: E \to E$ such that

$$\lim_{n \to \infty} \|F_n(x) - F(x)\|_{-\alpha} = 0$$

for all $x \in E$. For every $C > 0$, there exists a $K > 0$ such that for all $x \in E$ with $\|x\|_E \leq C$ and all $n \in \mathbb{N}$, we have $\|F_n(x)\|_{-\alpha} \leq K$. Furthermore, there is a $\gamma > 0$ such that the dissipativity bound

(2.7) $$\langle x^*, F_n(x+y) \rangle \leq -\gamma \|x\|_E$$

holds for every $x^* \in \partial \|x\|_E$ and every $x, y \in E$ with $\|x\|_E \geq C(1 + \|y\|_E)^N$.

As in [5], Example D.3, one can check that in the case $E = \mathcal{C}([0,1], \mathbb{R}^d)$, the elements of $\partial \|x\|_E$ can be characterized as follows: $x^* \in \partial \|x\|_E$ if and only if there exists a probability measure $|x^*|$ on $[0,1]$ with $\operatorname{supp} |x^*| \subseteq \{u \in [0,1] | |x(u)| = \|x\|_\infty\}$ and such that

(2.8) $$x^*(y) = \int \left\langle y(u), \frac{x(u)}{|x(u)|} \right\rangle |x^*|(du)$$

for every $y \in E$. Loosely speaking, the dissipativity condition in (A4) then states that the drift $F_n$ points inward for all locations $u \in [0,1]$, where $|x(u)|$ is largest and thus acts to decrease $\|x\|_E$.

DEFINITION 2.4. An $E$-valued and $(\mathcal{F}_t)$-adapted process $x$ is called a *mild solution* of equation (2.1) if almost surely

$$(2.9) \qquad x(t) = S(t)x_0 + \int_0^t S(t-s)F(x(s))\,ds + z(t) \qquad \forall t \geq 0$$

holds, where $z$ is the solution of the linear equation from (2.4).

LEMMA 2.5. *Let $\mathcal{L}$ satisfy assumptions* (A1) *and* (A2). *Let $F: E \to E^*$ be Lipschitz continuous on bounded sets, $\psi: \mathbb{R}_+ \to E$ be a continuous function and $x_0 \in \mathcal{H}^{1/2}$. Then, the equation*

$$(2.10) \qquad \frac{dx}{dt}(t) = \mathcal{L}x(t) + F(x(t) + \psi(t)), \qquad x(0) = x_0$$

*has a unique, local, $\mathcal{H}^{1/2}$-valued mild solution.*

*Furthermore, the length of the solution interval is bounded from below uniformly in $\|x_0\|_{1/2} + \sup_{t \in [0,1]} \|\psi(t)\|_E$.*

PROOF. Since $\psi$ is continuous, $\|\psi(t)\|_E$ is locally bounded. It is a straightforward exercise using (2.3) to show that, for sufficiently small $T$, the map $M_T$ acting from $\mathcal{C}([0,T], \mathcal{H}^{1/2})$ into itself and defined by

$$(M_T y)(t) = S(t)x_0 + \int_0^t S(t-s)F(y(s) + \psi(s))\,ds$$

is a contraction on a small ball around the element $t \mapsto S(t)x_0$. Therefore, (2.10) has a unique local solution in $\mathcal{H}^{1/2}$. The claim on the possible choice of $T$ can be checked in an equally straightforward way. □

THEOREM 2.6. *Let $\mathcal{L}$ and $F$ satisfy assumptions* (A1)–(A4). *Then, for every $x_0 \in E$, the equation (2.1) has a global, $E$-valued, unique mild solution and there exist positive constants $K_p$ and $\sigma$ such that*

$$\mathbb{E}\|x(t)\|_E^p \leq e^{-p\sigma t}\|x_0\|_E^p + K_p$$

*for all times.*

PROOF. Let $z$ be the solution of the linear equation $dz = \mathcal{L}z(t)\,dt + \sqrt{2}\,dw$ and, for $n \in \mathbb{N}$, let $y_n$ be the solution of

$$\frac{dy_n}{dt}(t) = \mathcal{L}y_n(t) + F_n(y_n(t) + z(t)), \qquad y_n(0) = x_0,$$

where $F_n$ is the approximation of $F$ from (A4). From Lemmas 2.1 and 2.5, we get that the differential equation almost surely has a local mild solution. We begin the proof by showing that $y_n$ can be extended to a global solution and obtaining an a priori bound for $y_n$ which does not depend on $n$.

Let $t \geq 0$ be sufficiently small that $y_n(t+h)$ exists for some $h > 0$. As an abbreviation, define $f(s) = F_n(y_n(s) + z(s))$ for all $s < t + h$. We then have

$$\|y(t+h)\|_E = \left\|S(h)y(t) + \int_t^{t+h} S(t+h-s)f(s)\,ds\right\|_E.$$

Since $f$ is continuous and the semigroup $S$ is a strongly continuous contraction semigroup on $E$, we obtain

$$\left\|\int_t^{t+h} S(t+h-s)f(s)\,ds - hS(h)f(t)\right\|_E$$

$$\leq \int_t^{t+h} \|S(t+h-s)(f(s)-f(t))\|_E + \|(S(t+h-s) - S(h))f(t)\|_E\,ds$$

$$\leq \int_t^{t+h} \|f(s) - f(t)\|_E\,ds + \int_0^h \|(S(r) - S(0))f(t)\|_E\,dr$$

$$= o(h)$$

and thus

$$\|y(t+h)\|_E = \|S(h)y(t) + S(h)hf(t)\|_E + o(h) \leq \|y(t) + hf(t)\|_E + o(h)$$

as $h \downarrow 0$. This gives

$$\limsup_{h \downarrow 0} \frac{\|y(t+h)\|_E - \|y(t)\|_E}{h}$$

$$\leq \lim_{h \downarrow 0} \frac{\|y(t) + hf(t)\|_E - \|y(t)\|_E}{h} = \max\{\langle y^*, f(t)\rangle | y^* \in \partial \|y(t)\|_E\},$$

where the last equation comes from [5], equation (D.2). Using assumption (A4), we get

$$\limsup_{h \downarrow 0} \frac{\|y_n(t+h)\|_E - \|y_n(t)\|_E}{h} \leq -\gamma \|y_n(t)\|_E$$

for all $t > 0$ with $\|y_n(t)\|_E \geq C(1 + \|z(t)\|_E)^N$.

An elementary proof shows that any continuous function $f : [0, T] \to \mathbb{R}$ with $f(t) > f(0)\exp(-\gamma t)$ for a $t \in (0, T]$ satisfies $\limsup(f(s+h) - f(s))/h > -\gamma f(s)$ for some $s \in [0, t)$. Therefore, whenever $\|y_n(t)\|_E \geq C(1 + \|z(t)\|_E)^N$ for all $t \in [a, b]$, the solution $y_n$ decays exponentially on this interval with

$$\|y_n(t)\|_E \leq \|y_n(a)\|_E e^{-\gamma(t-a)}$$

for all $t \in [a, b]$. Thus (see Figure 1), we find the a priori bound

(2.11)        $$\|y_n(t)\|_E \leq e^{-\gamma t}\|x_0\| \vee \sup_{0 < s < t} Ce^{-\gamma(t-s)}(1 + \|z(s)\|_E)^N$$

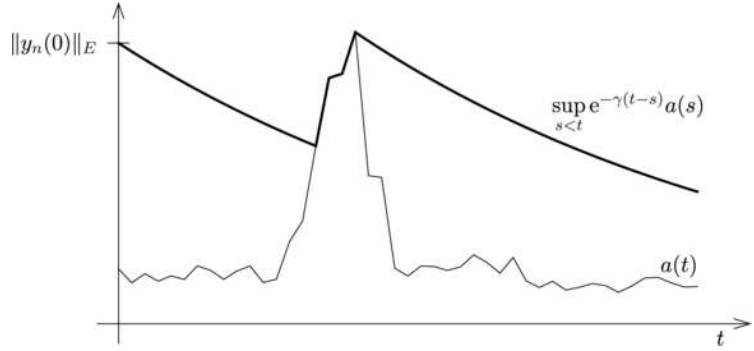

Fig. 1. *This illustrates the a priori bound on $\|y_n\|_E$ obtained in the proof of Theorem 2.6. Whenever $\|y_n(t)\|_E$ is above $a(t) = C(1 + \|z(t)\|_E)^N$, it decays exponentially. Therefore, the thick line is an upper bound for $\|y_n\|_E$.*

for the solution $y_n$. Using this bound and Lemma 2.5 repeatedly allows us to extend the solution $y_n$ to arbitrarily long time intervals.

Lemma 2.5 also gives local existence for the solution $y$ of

$$(2.12) \qquad \frac{dy}{dt}(t) = \mathcal{L}y(t) + F(y(t) + z(t)), \qquad y(0) = x_0.$$

Once we have seen that the bound (2.11) also holds for $y$, we obtain the required global existence for $y$. Let $t$ be sufficiently small for $y(t)$ to exist. Then, using (2.3),

$$\|y_n(s) - y(s)\|_E \leq C \left\| \int_0^s S(s-r)(F_n(y_n + z) - F(y+z))\, dr \right\|_\alpha$$
$$\leq C \int_0^s (s-r)^{-2\alpha} \|F_n(y_n + z) - F(y+z)\|_{-\alpha}\, dr$$

for every $s \leq t$ and thus

$$\int_0^t \|y_n - y\|_E\, ds \leq C \int_0^t \int_0^s (s-r)^{-2\alpha} \|F_n(y_n + z) - F(y_n + z)\|_{-\alpha}\, dr\, ds$$
$$+ C \int_0^t \int_0^s (s-r)^{-2\alpha} \|F(y_n + z) - F(y + z)\|_{-\alpha}\, dr\, ds$$
$$=: C(I_1 + I_2).$$

The map $F: E \to E^*$ is Lipschitz on bounded sets and thus has the same property when considered as a map $E \to \mathcal{H}^{-\alpha}$. Using (2.11) to see that there is a ball in $E$ which contains all $y_n$, we get $\|F(y_n + z) - F(y + z)\|_{-\alpha} \leq C\|y_n - y\|_{-\alpha}$. Fubini's theorem then gives

$$I_2 = \int_0^t \|y_n(r) - y(r)\|_{-\alpha} \int_r^t (s-r)^{-2\alpha}\, ds\, dr$$

$$\leq \frac{t^{1-2\alpha}}{1-2\alpha} \int_0^t \|y_n(r) - y(r)\|_E \, dr$$

and, by choosing $t$ sufficiently small and moving the $I_2$-term to the left-hand side, we find

$$\int_0^t \|y_n(r) - y(r)\|_E \, dr \leq CI_1.$$

By (A4), the term $\|F_n(y_n + z) - F(y_n + z)\|_{-\alpha}$ in the integral is uniformly bounded by some constant $K$ and thus $(s-r)^{-2\alpha}M$ is an integrable, uniform upper bound for the integrand. Again by (A4), the integrand converges to 0 pointwise, so the dominated convergence theorem yields

(2.13) $$\int_0^t \|y_n(r) - y(r)\|_E \, dr \leq CI_1 \longrightarrow 0$$

as $n \to \infty$. Assume (for the purposes of obtaining a contradiction) that $y$ violates the bound (2.11) for some time $s \in [0, t]$. Since $t \mapsto y(t)$ is continuous, the bound is violated for a time interval of positive length, so $\int_0^t \|y_n(r) - y(r)\|_E \, dr$ is bounded from below uniformly in $n$. This contradicts (2.13), so $y$ must satisfy the a priori estimate (2.11). Again, we can iterate this step and extend the solution $y$ of (2.12) and thus the solution $x = y + z$ of (2.1) to arbitrary large times.

Now, all that remains, is to prove the given bound on $\mathbb{E}\|x(t)\|_E^p$. For $k \in \mathbb{N}$, let $a_k = \sup_{k-1 \leq t \leq k} \|z(t)\|_E$ and

$$\xi_k = \sup_{k+1 \leq t \leq k+2} \sqrt{2} \left\| \int_k^t S(s-k) \, dw(s) \right\|_E.$$

The $\xi_k$ are then identically distributed and, for $|k - l| \geq 2$, the random variables $\xi_l$ and $\xi_k$ are independent. Without loss of generality, we can assume $\|S(t)x\|_E \leq e^{-t\varepsilon}\|x\|_E$ for some small value $\varepsilon > 0$ [otherwise, replace $\mathcal{L}$ with $\mathcal{L} - \varepsilon I$ and $F$ with $F + \varepsilon I$, where $\varepsilon$ is chosen sufficiently small that (A4) still holds for $F + \varepsilon I$]. Thus, for $h \in [1, 2]$, we get

$$\|z(k+h)\| \leq \|S(h)z(k)\|_E + \sqrt{2} \left\| \int_k^{k+h} S(s-k) \, dw(s) \right\|_E \leq e^{-\varepsilon} a_k + \xi_k$$

and, consequently, $a_{k+2} \leq e^{-\varepsilon} a_k + \xi_k$. Since the $\xi_k$ and $a_1$, $a_2$ have Gaussian tails, it is a straightforward calculation to check from this recursion relation that the expression $\sum_{k=1,\ldots,m} e^{\gamma(m-k)} a_k^N$ has bounded moments of all orders that are independent of $m$. Since the right-hand side of (2.11) is bounded by expressions of this type, the required bound on the solutions $x(t)$ follows immediately, with $\sigma = \gamma - \varepsilon$. □

2.2. *The preconditioned case.* In this section, we consider semilinear SPDEs of the form

$$(2.14) \qquad dx = \mathcal{G}(\mathcal{L}x + F(x))\,dt + \sqrt{2}\mathcal{G}^{1/2}\,dw(t), \qquad x(0) = x_0,$$

where $\mathcal{L}$, $F$ and $w$ are as before and $\mathcal{G}$ is a self-adjoint, positive linear operator on $\mathcal{H}$. We seek a strong solution of this equation, defined below. In order to simplify our notation, we define $\tilde{\mathcal{L}} = \mathcal{G}\mathcal{L}$, $\tilde{F} = \mathcal{G}F$ and $\tilde{w} = \mathcal{G}^{1/2}w$. Then, $\tilde{w}$ is a $\mathcal{G}$-Wiener process on $\mathcal{H}$ and equation (2.14) can be written as

$$dx = \tilde{\mathcal{L}}x\,dt + \tilde{F}(x)\,dt + \sqrt{2}\,d\tilde{w}(t), \qquad x(0) = x_0.$$

For the operator $\mathcal{L}$, we will continue to use assumptions (A1) and (A2). For $F$, we use the growth condition (A3), but replace the dissipativity condition (A4) with the following one.

(A5) There exists $N > 0$ such that $F$ satisfies

$$\langle x, F(x+y)\rangle \le C(1 + \|y\|_E)^N$$

for every $x, y \in E$.

REMARK 2.7. Note that (A5) is structurally similar to assumption (A4) above, except that we now assume dissipativity in $\mathcal{H}$ rather than in $E$.

We make the following assumption on $\mathcal{G}$:

(A6) The operator $\mathcal{G}: \mathcal{H} \to \mathcal{H}$ is trace class, self-adjoint and positive definite, the range of $\mathcal{G}$ is dense in $\mathcal{H}$ and the Gaussian measure $\mathcal{N}(0, \mathcal{G})$ is concentrated on $E$.

Define the space $\tilde{\mathcal{H}}$ to be $\mathcal{D}(\mathcal{G}^{-1/2})$ with the inner product $\langle x, y\rangle_{\tilde{\mathcal{H}}} = \langle x, \mathcal{G}^{-1}y\rangle$. We then assume that $\mathcal{G}$ is equal to the inverse of $\mathcal{L}$, up a "small" error in the following sense.

(A7) We have $\mathcal{G}\mathcal{L} = -I + \mathcal{K}$, where $\mathcal{K}$ is a bounded operator from $\mathcal{H}$ to $\tilde{\mathcal{H}}$.

LEMMA 2.8. *Assume* (A1), (A2), (A6), (A7). *Then, $\tilde{\mathcal{H}} = \mathcal{H}^{1/2}$. In particular, $\tilde{\mathcal{H}} \subset E$.*

PROOF. First, note that by [24], Theorem VII.1.3, the fact that $\mathcal{G}\mathcal{L}$ is bounded on $\mathcal{H}$ implies that $\text{range}(\mathcal{G}) \subset \mathcal{D}(\mathcal{L})$. Furthermore, (A7) implies that $\mathcal{D}(\mathcal{L}) \subset \tilde{\mathcal{H}}$. For every $x \in \text{range}\,\mathcal{G}$, one has

$$|\|x\|_{\tilde{\mathcal{H}}}^2 - \|x\|_{1/2}^2| = |\langle \mathcal{G}^{-1/2}x, \mathcal{G}^{-1/2}\mathcal{K}x\rangle|$$
$$\le \|x\|_{\tilde{\mathcal{H}}}\|\mathcal{K}x\|_{\tilde{\mathcal{H}}} \le C\|x\|_{\tilde{\mathcal{H}}}\|x\| \le C\|x\|_{\tilde{\mathcal{H}}}\|x\|_{1/2},$$

so the norms $\|x\|_{\tilde{\mathcal{H}}}$ and $\|x\|_{1/2}$ are equivalent. In particular, we have

$$\text{range}(\mathcal{G}) \subset \mathcal{D}(\mathcal{L}) \subset \tilde{\mathcal{H}} \subset \mathcal{H}^{1/2}.$$

The facts that range$(\mathcal{G})$ is dense in $\tilde{\mathcal{H}}$ and $\mathcal{D}(\mathcal{L})$ is dense in $\mathcal{H}^{1/2}$ conclude the proof. $\square$

DEFINITION 2.9. An $E$-continuous and adapted process $x$ is called a *strong solution* of (2.14) if it satisfies

$$(2.15) \quad x(t) = x_0 + \int_0^t (\mathcal{G}\mathcal{L}x(s) + \mathcal{G}F(x(s)))\,ds + \sqrt{2}\tilde{w}(t) \qquad \forall t \geq 0$$

almost surely.

THEOREM 2.10. *Let $\tilde{\mathcal{L}}$, $\tilde{F}$ and $\mathcal{G}$ satisfy assumptions* (A1)–(A3) *and* (A5)–(A7). *Then, for every $x_0 \in E$, equation* (2.14) *has a global, $E$-valued, unique strong solution. There exists a constant $N > 0$ and, for every $p > 0$, there exist constants $K_p$, $C_p$ and $\gamma_p > 0$ such that*

$$(2.16) \qquad \mathbb{E}\|x(t)\|_E^p \leq C_p(1 + \|x_0\|_E)^{Np} e^{-\gamma_p t} + K_p$$

*for all times.*

PROOF. Since it follows from (A6) and Kolmogorov's continuity criterion that the process $\tilde{w}(t)$ is $E$-valued and has continuous sample paths, it is a straightforward exercise (use Picard iterations pathwise) to show that (2.14) has a unique strong solution lying in $E$ for all times. It is possible to obtain uniform bounds on this solution in the following way. Choose an arbitrary initial condition $x_0 \in E$ and let $y$ be the solution to the linear equation

$$dy = -y\,dt + d\tilde{w}(t), \qquad y(0) = x_0.$$

There exist constants $\tilde{K}_p$ such that

$$(2.17) \qquad \mathbb{E}\|y(t)\|_E^p \leq e^{-pt}\|x_0\|_E^p + \tilde{K}_p.$$

Denote by $z$ the difference $z(t) = x(t) - y(t)$. It then follows that $z$ satisfies the ordinary differential equation

$$\frac{dz}{dt} = \tilde{\mathcal{L}}z(t) + \tilde{F}(x(t)) + \mathcal{K}y(t), \qquad z(0) = 0.$$

Since $\tilde{\mathcal{L}}$ is bounded from $\tilde{\mathcal{H}}$ to $\tilde{\mathcal{H}}$ by (A7) and $\tilde{F}(x) + \mathcal{K}y \in \tilde{\mathcal{H}}$ for every $x, y \in E$ by Lemma 2.8, it follows that $z(t) \in \tilde{\mathcal{H}}$ for all times. Furthermore,

SPDEs ARISING IN PATH SAMPLING 15we have the following bound on its moments:

$$\frac{d\|z\|_{\tilde{\mathcal{H}}}^2}{dt} \leq -2w\|z\|_{\tilde{\mathcal{H}}}^2 + \langle \tilde{F}(x), x - y \rangle_{\tilde{\mathcal{H}}} + \langle \mathcal{K}y, z \rangle_{\tilde{\mathcal{H}}}$$

$$\leq C - 2w\|z\|_{\tilde{\mathcal{H}}}^2 + \frac{w}{2}\|z\|_{\tilde{\mathcal{H}}}^2 + C(1+\|y\|_E)^N + \frac{1}{2w}\|\mathcal{K}y\|_{\tilde{\mathcal{H}}}^2$$

$$\leq -w\|z\|_{\tilde{\mathcal{H}}}^2 + C(1+\|y\|_E)^N.$$

Using Gronwall's lemma, it thus follows from (2.17) that $x$ satisfies a bound of the type (2.16) for every $p \geq 0$. $\square$

**3. Stationary distributions of semilinear SPDEs.** In this section, we give an explicit representation of the stationary distribution of (2.1) and (2.14) when $F$ is a gradient, by comparing it to the stationary distribution of the linear equation

(3.1)
$$\frac{dz}{dt}(t) = \mathcal{L}z(t) + \sqrt{2}\frac{dw}{dt}(t) \qquad \forall t \geq 0,$$
$$z(0) = 0.$$

The main results are stated in Theorems 3.4 and Theorem 3.6.

The solution of (3.1) is the process $z$ from Lemma 2.1 and its stationary distribution is the Gaussian measure $\nu = \mathcal{N}(0, -\mathcal{L}^{-1})$. In this section, we identify, under the assumptions of Section 2 and with $F = U'$ for a Fréchet differentiable function $U : E \to \mathbb{R}$, the stationary distribution of the equations (2.1) and (2.14). It transpires to be the measure $\mu$ which has the Radon–Nikodym derivative

$$d\mu = c\exp(U)\,d\nu$$

with respect to the stationary distribution $\nu$ of the linear equation, where $c$ is the appropriate normalization constant. In the next section, we will see that there are no other stationary distributions.

The results here are slight generalizations of the results in [25]. Our situation differs from the one in [25] in that we allow the nonlinearity $U'$ to take values in $E^*$ instead of $\mathcal{H}$ and that we consider preconditioning for the SPDE. We have scaled the noise by $\sqrt{2}$ to simplify notation. Where possible, we refer to the proofs in [25] and describe in detail arguments which require nontrivial extensions of that paper.

Let $(e_n)_{n \in \mathbb{N}}$, be an orthonormal set of eigenvectors of $\mathcal{L}$ in $\mathcal{H}$. For $n \in \mathbb{N}$ let $E_n$ be the subspace spanned by $e_1, \ldots, e_n$ and let $\Pi_n$ be the orthogonal projection onto $E_n$. From [25], Proposition 2, we know that, under assumption (A2), we have $E_n \subseteq E$ for every $n \in \mathbb{N}$.

LEMMA 3.1. *Suppose that assumptions* (A1) *and* (A2) *are satisfied. There then exist linear operators* $\hat{\Pi}_n \colon E \to E_n$ *which are uniformly bounded in the operator norm on $E$ and which satisfy* $\hat{\Pi}_n \Pi_n = \hat{\Pi}_n$ *and* $\|\hat{\Pi}_n x - x\|_E \to 0$ *as $n \to \infty$.*

PROOF. The semigroup $S$ on $\mathcal{H}$ can be written as

$$S(t)x = \sum_{k=1}^{\infty} e^{-t\lambda_k} \langle e_k, x \rangle e_k$$

for all $x \in \mathcal{H}$ and $t \geq 0$ where the series converges in $\mathcal{H}$. Since there is a constant $c_1 > 0$ with $\|x\|_{\mathcal{H}} \leq c_1 \|x\|_E$ and, from [25], Proposition 2, we know there exists a constant $c_2 > 0$ with $\|e_k\|_E \leq c_2 \sqrt{\lambda_k}$, we have

$$\|e^{-t\lambda_k} \langle e_k, x \rangle e_k\|_E \leq e^{-t\lambda_k} \|e_k\|_{\mathcal{H}} \|x\|_{\mathcal{H}} \|e_k\|_E \leq c_1 c_2 e^{-t\lambda_k} \sqrt{\lambda_k} \|x\|_E$$

for every $k \in \mathbb{N}$. Consequently, there is a constant $c_3 > 0$ with

$$\|e^{-t\lambda_k} \langle e_k, x \rangle e_k\|_E \leq c_3 t^{-3/2} \lambda_k^{-1} \|x\|_E.$$

Now, define $\hat{\Pi}_n$ by

(3.2) $$\hat{\Pi}_n x = \sum_{k=1}^{n} e^{-t_n \lambda_k} \langle e_k, x \rangle e_k,$$

where

$$t_n = \left( \sum_{k=n+1}^{\infty} \lambda_k^{-1} \right)^{1/3}.$$

[This series converges, since assumption (A2) implies that $\mathcal{L}^{-1}$ is trace class.] Then,

$$\|(S(t_n) - \hat{\Pi}_n)x\|_E \leq c_3 t_n^{-3/2} \sum_{k=n+1}^{\infty} \lambda_k^{-1} \|x\|_E = c_3 t_n^{3/2} \|x\|_E.$$

We have $\|\hat{\Pi}_n\|_E \leq \|S(t_n) - \hat{\Pi}_n\|_E + \|S(t_n)\|_E$. Since $S$ is strongly continuous on $E$, the norms $\|S(t_n)\|_E$ are uniformly bounded. Thus, the operators $\hat{\Pi}_n$ are uniformly bounded and, since $t_n \to 0$, we have $\|\hat{\Pi}_n x - x\|_E \leq \|S(t_n)x - \hat{\Pi}_n x\|_E + \|S(t_n)x - x\|_E \to 0$ as $n \to \infty$. □

Since the eigenvectors $e_n$ are contained in each of the spaces $\mathcal{H}^\alpha$, we can consider $\hat{\Pi}_n$, as defined by (3.2), to be an operator between any two of the spaces $E$, $E^*$, $\mathcal{H}$ and $\mathcal{H}^\alpha$ for all $\alpha \in \mathbb{R}$. In the sequel, we will simply write $\hat{\Pi}_n$ for all of these operators. Taken from $\mathcal{H}$ to $\mathcal{H}$, this operator is self-adjoint. The adjoint of the operator $\hat{\Pi}_n$ from $E$ to $E$ is just the $\hat{\Pi}_n$ we obtain by using (3.2) to define an operator from $E^*$ to $E^*$. Therefore, in our notation, we never need to write $\hat{\Pi}_N^*$. As a consequence of Lemma 3.1, the operators $\hat{\Pi}_n$ are uniformly bounded from $E^*$ to $E^*$.

Denote the space of bounded, continuous functions from $E$ to $\mathbb{R}$ by $\mathcal{C}_b(E)$. We state and prove a modified version of [25], Theorem 2.

THEOREM 3.2. *Suppose that assumptions* (A1), (A2) *are satisfied. Let $\mathcal{G}$ be a positive definite, self-adjoint operator on $\mathcal{H}$, let $U : E \to \mathbb{R}$ be bounded from above and Fréchet-differentiable and, for $n \in \mathbb{N}$, let $(P_t^n)_{t>0}$ be the semigroup on $\mathcal{C}_b(E)$ which is generated by the solutions of*

$$(3.3) \qquad dx(t) = \mathcal{G}_n(\mathcal{L}x + F_n(x(t)))\, dt + \sqrt{2}\mathcal{G}_n^{1/2}\Pi_n\, dw,$$

*where $U_n = U \circ \hat{\Pi}_n$, $F_n = U_n'$, $\mathcal{G}_n = \hat{\Pi}_n \mathcal{G} \hat{\Pi}_n$ and $w$ is a cylindrical Wiener process. Define the measure $\mu$ by*

$$d\mu(x) = e^{U(x)}\, d\nu(x),$$

*where $\nu = \mathcal{N}(0, -\mathcal{L}^{-1})$. Let $(P_t)_{t>0}$ be a semigroup on $\mathcal{C}_b(E)$ such that $P_t^n \varphi(x_n) \to P_t \varphi(x)$ for every $\varphi \in \mathcal{C}_b(E)$, for every sequence $(x_n)$ with $x_n \in E_n$ and $x_n \to x \in E$ and for every $t > 0$. The semigroup $(P_t)_{t>0}$ is then $\mu$-symmetric.*

PROOF. From [25], Theorem 1, we know that the stationary distribution of $z$ is $\nu$ and, from the finite dimensional theory, we know that (3.3) is reversible with a stationary distribution $\mu_n$ which is given by

$$d\mu_n(x) = c_n e^{U_n(x)}\, d\nu_n(x),$$

where $\nu_n = \nu \circ \Pi_n^{-1}$ and $c_n$ is the appropriate normalization constant. Thus, for all continuous, bounded $\varphi, \psi : E \to \mathbb{R}$, we have

$$\int_E \varphi(x) P_t^n \psi(x)\, d\mu_n(x) = \int_E \psi(x) P_t^n \varphi(x)\, d\mu_n(x)$$

and substitution gives

$$\int_E \varphi(\Pi_n x) P_t^n \psi(\Pi_n x) e^{U(\hat{\Pi}_n x)}\, d\nu(x)$$

$$= \int_E \psi(\Pi_n x) P_t^n \varphi(\Pi_n x) e^{U(\hat{\Pi}_n x)}\, d\nu(x)$$

for every $t \geq 0$ and every $n \in \mathbb{N}$.

As in the proof of [25], Theorem 2, we obtain $\Pi_n x \to x$ in $E$ for $\nu$-a.a. $x$. Since $U$ is bounded from above and continuous and $\varphi, \psi \in \mathcal{C}_b(E)$, we can use the dominated convergence theorem to conclude

$$\int_E \varphi(x) P_t \psi(x) e^{U(x)}\, d\nu(x) = \int_E \psi(x) P_t \varphi(x) e^{U(x)}\, d\nu(x).$$

This shows that the semigroup $(P_t)_{t>0}$ is $\mu$-symmetric. □

3.1. *The nonpreconditioned case.* We will apply Theorem 3.2 in two different situations, namely for $\mathcal{G} = I$ (in this subsection) and for $\mathcal{G} \approx -\mathcal{L}^{-1}$ (in the next subsection). The case $\mathcal{G} = I$ is treated in [25], Proposition 5 and [25], Theorem 4. Since, in the present text, we allow the nonlinearity $U'$ to take values in $E^*$ instead of $\mathcal{H}$, we repeat the (slightly modified) result here.

LEMMA 3.3. *For $n \in \mathbb{N}$, let $F_n, F : E \to E^*$, $T > 0$ and $\psi_n, \psi : [0, T] \to E$ be continuous functions such that the following conditions hold:*

- *for every $r > 0$, there exists a $K_r > 0$ such that $\|F_n(x) - F_n(y)\|_{E^*} \le K_r \|x - y\|_E$ for every $x, y \in E$ with $\|x\|_E, \|y\|_E \le r$ and every $n \in \mathbb{N}$;*
- $F_n(x) \to F(x)$ *in $E^*$ as $n \to \infty$ for every $x \in E$;*
- $\psi_n \to \psi$ *in $\mathcal{C}([0, T], E)$ as $n \to \infty$;*
- *there exists a $p > 1$ with*

$$\text{(3.4)} \qquad \int_0^T \|S(s)\|_{E^* \to E}^p \, ds < \infty.$$

*Let $u_n, u : [0, T] \to E$ be the solutions of*

$$\text{(3.5)} \qquad u_n(t) = \int_0^t S(t - s) F_n(u_n(s)) \, ds + \psi_n(t),$$

$$\text{(3.6)} \qquad u(t) = \int_0^t S(t - s) F(u(s)) \, ds + \psi(t).$$

*Then, $u_n \to u$ in $\mathcal{C}([0, T], E)$.*

PROOF. We have

$$\|u_n(t) - u(t)\|_E \le \left\| \int_0^t S(t - s)(F_n(u(s)) - F(u(s))) \, ds \right\|_E$$
$$+ \left\| \int_0^t S(t - s)(F_n(u_n(s)) - F_n(u(s))) \, ds \right\|_E$$
$$+ \|\psi_n(t) - \psi(t)\|_E$$
$$= I_1(t) + I_2(t) + I_3(t)$$

for all $t \in [0, T]$. We can choose $q > 1$ with $1/p + 1/q = 1$ to obtain

$$I_1(t) \le \int_0^t \|S(t - s)(F_n(u(s)) - F(u(s)))\|_E \, ds$$
$$\le \int_0^t \|S(t - s)\|_{E^* \to E} \|F_n(u(s)) - F(u(s))\|_{E^*} \, ds$$
$$\le \left( \int_0^T \|S(t - s)\|_{E^* \to E}^p \, ds \right)^{1/p} \left( \int_0^T \|F_n(u(s)) - F(u(s))\|_{E^*}^q \, ds \right)^{1/q}.$$

By dominated convergence, the right-hand side converges to 0 uniformly in $t$ as $n \to \infty$.

For $n \in \mathbb{N}$ and $r > 0$, define
$$\tau_{n,r} = \inf\{t \in [0,T] | \|u(t)\|_E \geq r \text{ or } \|u_n(t)\|_E \geq r\},$$
with the convention that $\inf \varnothing = T$. For $t \leq \tau_{n,r}$ we have
$$I_2(t) \leq K_r \int_0^t \|S(t-s)\|_{E^* \to E} \|u_n(s) - u(s)\|_E \, ds$$
and, consequently,
$$\|u_n(t) - u(t)\|_E \leq \sup_{0 \leq t \leq T} I_1(t) + \|\psi_n(t) - \psi(t)\|_E$$
$$+ K_r \int_0^t \|S(t-s)\|_{E^* \to E} \|u_n(s) - u(s)\|_E \, ds.$$

Using Gronwall's lemma, we can conclude that
$$\|u_n(t) - u(t)\|_E \leq \left(\sup_{0 \leq t \leq T} I_1(t) + \|\psi_n(t) - \psi(t)\|_E\right)$$
$$\times \exp\left(K_r \int_0^T \|S(s)\|_{E^* \to E} \, ds\right)$$
for all $t \leq \tau_{n,r}$.

Now, choose $r > 0$ such that $\sup_{0 \leq t \leq T} \|u(t)\|_E \leq r/2$. Then, for sufficiently large $n$ and all $t \leq \tau_{n,r}$, we have $\|u_n(t) - u(t)\|_E \leq r/2$ and thus $\sup_{0 \leq t \leq T} \|u(t)\|_E \leq r$. This implies that $\tau_{n,r} = T$ for sufficiently large $n$ and the result follows. □

With all of these preparations in place, we can now show that the measure $\mu$ is a stationary distribution of the nonpreconditioned equation. The proof works by approximating the infinite-dimensional solution of (2.1) by finite-dimensional processes. Lemma 3.3 then shows that these finite dimensional processes converge to the solution of (2.1) and Theorem 3.2 finally shows that the corresponding stationary distributions also converge.

THEOREM 3.4. *Let $U : E \to \mathbb{R}$ be bounded from above and Fréchet differentiable. Assume that $\mathcal{L}$ and $F = U'$ satisfy assumptions* (A1)–(A4). *Define the measure $\mu$ by*

(3.7) $$d\mu(x) = ce^{U(x)} \, d\nu(x),$$

*where $\nu = \mathcal{N}(0, -\mathcal{L}^{-1})$ and $c$ is a normalization constant. Then,* (2.1) *has a unique mild solution for every initial condition $x_0 \in E$ and the corresponding Markov semigroup on $E$ is $\mu$-symmetric. In particular, $\mu$ is an invariant measure for* (2.1).

PROOF. Let $x_0 \in E$. From Theorem 2.6, the SDE (2.1) has a mild solution $x$ starting at $x_0$. Defining $\psi(t) = S(t)x_0 + z(t)$, where $z$ is given by (3.1), we can a.s. write this solution in the form (3.6). Now, consider a sequence $(x_0^n)$ with $x_0^n \in E_n$ for all $n \in \mathbb{N}$ and $x_0^n \to x_0$ as $n \to \infty$. Let $\mathcal{G} = I$. Then, for every $n \in \mathbb{N}$, the finite-dimensional equation (3.3) has a solution $x^n$ which starts at $x_0^n$ and this solution can a.s. be written in the form (3.5), with $\psi_n = S(t)x_0^n + z_n(t)$ and $z_n = \Pi_n z$. From [25], Proposition 1, we get that $z_n \to z$ as $n \to \infty$ and thus $\psi_n \to \psi$ in $\mathcal{C}([0,T], E)$ as $n \to \infty$.

Define $F_n$ as in Theorem 3.2. We then have $F_n(x) = \hat{\Pi}_n F(\hat{\Pi}_n x)$ and thus $F_n(x) \to F(x)$ as $n \to \infty$ for every $x \in E$. Also, since $F$ is locally Lipschitz, and $\hat{\Pi}_n : E \to E$ and $\hat{\Pi}_n : E^* \to E^*$ are uniformly bounded, the $F_n$ are locally Lipschitz, where the constant can be chosen uniformly in $n$. From (2.3), we obtain
$$\int_0^T \|S(t)\|_{E^* \to E} \, dt \leq c_2 \int_0^T t^{-2\alpha} \, dt < \infty$$
for every $T > 0$ and thus condition (3.4) is satisfied. We can now use Lemma 3.3 to conclude that $x^n \to x$ in $\mathcal{C}([0,T], E)$ as $n \to \infty$ almost surely. Using dominated convergence, we see that $P_t^n \varphi(x_n) \to P_t \varphi(x)$ for every $\varphi \in \mathcal{C}_b(E)$ and every $t > 0$, where $(P_t^n)$ are the semigroups from Theorem 3.2 and $(P_t)_{t>0}$ is the semigroup generated by the solutions of (2.1). We can now apply Theorem 3.2 to conclude that $(P_t)_{t>0}$ is $\mu$-symmetric. $\square$

3.2. *The preconditioned case.* For the preconditioned case, we require the covariance operator $\mathcal{G}$ of the noise to satisfy assumptions (A6) and (A7), in particular for $\mathcal{G}$ to be trace class. Thus, we can use strong solutions of (3.3) here. The analogue of Lemma 3.3 is given in the following lemma.

LEMMA 3.5. *Let $T > 0$ and, for $n \in \mathbb{N}$, let $\tilde{\mathcal{L}}_n, \tilde{\mathcal{L}}$ be bounded operators on $E$ and let $\tilde{F}_n, \tilde{F} : E \to E$ as well as $\psi_n, \psi : [0,T] \to E$ be continuous functions such that the following conditions hold:*

- $\tilde{\mathcal{L}}_n x \to \tilde{\mathcal{L}} x$ and $\tilde{F}_n(x) \to \tilde{F}(x)$ in $E$ as $n \to \infty$ for every $x \in E$;
- *for every $r > 0$, there is a $K_r > 0$ such that*

(3.8) $$\|\tilde{F}_n(x) - \tilde{F}_n(y)\|_E \leq K_r \|x - y\|_E$$

  *for every $x, y \in E$ with $\|x\|_E, \|y\|_E \leq r$ and every $n \in \mathbb{N}$;*
- $\psi_n \to \psi$ *in* $\mathcal{C}([0,T], E)$ *as* $n \to \infty$.

*Let $u_n, u : [0,T] \to E$ be solutions of*

(3.9) $$u_n(t) = \int_0^t (\tilde{\mathcal{L}}_n u_n(s) + \tilde{F}_n(u_n(s))) \, ds + \psi_n(t),$$

(3.10) $$u(t) = \int_0^t (\tilde{\mathcal{L}} u(s) + \tilde{F}(u(s))) \, ds + \psi(t),$$

*then $u_n \to u$ in $\mathcal{C}([0,T], E)$.*

PROOF. We have

$$\|u_n(t) - u(t)\|_E$$
$$\leq \int_0^t \|\tilde{\mathcal{L}}_n u(s) - \tilde{\mathcal{L}} u(s) + \tilde{F}_n(u(s)) - \tilde{F}(u(s))\|_E \, ds$$
$$+ \int_0^t \|\tilde{\mathcal{L}}_n u_n(s) - \tilde{\mathcal{L}}_n u(s) + \tilde{F}_n(u_n(s)) - \tilde{F}_n(u(s))\|_E \, ds$$
$$+ \|\psi_n(t) - \psi(t)\|_E$$
$$= I_1(t) + I_2(t) + I_3(t)$$

for all $t \in [0,T]$. By the uniform boundedness principle, we have $\sup_{n \in \mathbb{N}} \|\tilde{\mathcal{G}}_n\|_E < \infty$ and thus we can choose $K_r$ sufficiently large to obtain

$$\|\tilde{\mathcal{L}}_n x - \tilde{\mathcal{L}}_n y + \tilde{F}_n(x) - \tilde{F}_n(y)\|_E \leq K_r \|x - y\|_E$$

for every $x, y \in E$ with $\|x\|_E, \|y\|_E \leq r$ and every $n \in \mathbb{N}$. We also have

$$\sup_{0 \leq t \leq T} I_1(t) \leq \int_0^T \|\tilde{\mathcal{L}}_n u(s) - \tilde{\mathcal{L}} u(s) + \tilde{F}_n(u(s)) - \tilde{F}(u(s))\|_E \, ds \longrightarrow 0$$

as $n \to \infty$, by dominated convergence.

For $n \in \mathbb{N}$ and $r > 0$, define

$$\tau_{n,r} = \inf\{t \in [0,T] | \|u(t)\|_E \geq r \text{ or } \|u_n(t)\|_E \geq r\},$$

with the convention that $\inf \varnothing = T$. For $t \leq \tau_{n,r}$, we have

$$I_2(t) \leq K_r \int_0^t \|u_n(s) - u(s)\|_E \, ds$$

and, consequently,

$$\|u_n(t) - u(t)\|_E \leq \sup_{0 \leq t \leq T} I_1(t) + K_r \int_0^t \|u_n(s) - u(s)\|_E \, ds + \|\psi_n(t) - \psi(t)\|_E.$$

Using Gronwall's lemma, we can conclude that

$$\|u_n(t) - u(t)\|_E \leq e^{K_r T} \left( \sup_{0 \leq t \leq T} I_1(t) + \|\psi_n(t) - \psi(t)\|_E \right)$$

for all $t \leq \tau_{n,r}$.

Now, choose $r > 0$ such that $\sup_{0 \leq t \leq T} \|u(t)\|_E \leq r/2$. For sufficiently large $n$ and all $t \leq \tau_{n,r}$, we then have $\|u_n(t) - u(t)\|_E \leq r/2$ and thus $\sup_{0 \leq t \leq T} \|u(t)\|_E \leq r$. This implies that $\tau_{n,r} = T$ for sufficiently large $n$ and the result follows. $\square$

The following theorem shows that the measure $\mu$ is now also a stationary distribution of the preconditioned equation. Again the proof works by approximating the infinite-dimensional solution of (2.1) by finite-dimensional processes.

THEOREM 3.6. *Let $U: E \to \mathbb{R}$ be bounded from above and Fréchet differentiable. Assume that the operators $\mathcal{G}$ and $\mathcal{L}$ and the drift $F = U'$ satisfy assumptions* (A1)–(A3), *and* (A5)–(A7). *Define the measure $\mu$ by*

$$(3.11) \qquad d\mu(x) = ce^{U(x)}\,d\nu(x),$$

*where $\nu = \mathcal{N}(0, -\mathcal{L}^{-1})$ and $c$ is a normalization constant. Equation* (2.14) *then has a unique strong solution for every initial condition $x_0 \in E$ and the corresponding semigroup on $E$ is $\mu$-symmetric. In particular, $\mu$ is an invariant measure for* (2.14).

PROOF. Let $x_0 \in E$. From Theorem 2.10, SDE (2.14) has a strong solution $x$ starting at $x_0$. Defining $\psi(t) = x_0 + \tilde{w}(t)$, where $\tilde{w} = \mathcal{G}^{1/2}w$ is a $\mathcal{G}$-Wiener process, we can a.s. write this solution in the form (3.10). Now, consider a sequence $(x_0^n)$ with $x_0^n \in E_n$ for all $n \in \mathbb{N}$ and $x_0^n \to x_0$ as $n \to \infty$. For every $n \in \mathbb{N}$, the finite-dimensional equation (3.3) then has a solution $x^n$ which starts at $x_0^n$ and this solution can a.s. be written in the form (3.9), with $\psi_n = x_0^n + \hat{\Pi}_n \mathcal{G}^{1/2}w(t)$. Since the function $\psi$ is continuous, it can be approximated arbitrarily well by a piecewise affine function $\hat{\psi}$. Since the operators $\hat{\Pi}_n$ are equibounded in $E$ and satisfy $\hat{\Pi}_n y \to y$ for every $y \in E$, it is easy to see that $\hat{\Pi}_n \hat{\psi} \to \hat{\psi}$ in $\mathcal{C}([0,T], E)$. On the other hand, $\|\psi_n - \hat{\Pi}_n \hat{\psi}\|_E$ is bounded by $\|\hat{\Pi}_n x_0 - x_0^n\|_E + \|\hat{\Pi}_n\|_{E \to E} \|\psi - \hat{\psi}\|_E$, so it also gets arbitrarily small. This shows that $\psi_n$ indeed converges to $\psi$ in $\mathcal{C}([0,T], E)$.

Because of (A6) and (A7), we have $\|\mathcal{G}\|_{E \to E} < \infty$ and $\|\mathcal{G}\mathcal{L}\|_{E \to E} < \infty$. Let $F = U'$ and define $F_n$ and $\mathcal{G}_n$ as in Theorem 3.2. We then have $F_n(x) = \hat{\Pi}_n F(\hat{\Pi}_n x)$. Let $\tilde{\mathcal{L}}_n = \mathcal{G}_n \mathcal{L} = \hat{\Pi}_n \mathcal{G} \mathcal{L} \hat{\Pi}_n$, $\tilde{\mathcal{L}} = \mathcal{G}\mathcal{L}$, $\tilde{F}_n = \mathcal{G}_n F_n$ and $\tilde{F} = \mathcal{G}F$. Since $\|\hat{\Pi}_n\|_{E \to E} \le c$ for all $n \in \mathbb{N}$ and some constant $c < \infty$ and since $\|\hat{\Pi}_n x_n - x\|_E \le \|\hat{\Pi}_n x_n - \hat{\Pi}_n x\|_E + \|\hat{\Pi}_n x - x\|_E$, we have $\hat{\Pi}_n x_n \to x$ in $E$ as $n \to \infty$ for every sequence $(x_n)$ with $x_n \to x$ in $E$. Since $\mathcal{G}\mathcal{L}$ is a bounded operator on $E$, we can use this fact to obtain $\tilde{\mathcal{L}}_n x \to \tilde{\mathcal{L}} x$ in $E$ as $n \to \infty$ for every $x \in E$. Since $\mathcal{G}\mathcal{L}$ is bounded from $E$ to $E$ and $\mathcal{L}(E) \supseteq \mathcal{L}(\mathcal{H}^{1/2}) = \mathcal{H}^{-1/2}$, the operator $\mathcal{G}$ is defined on all of $E^* \subseteq \mathcal{H}^{-1/2}$ and thus bounded from $E^*$ to $E$ and we obtain $\tilde{F}_n(x) \to \tilde{F}(x)$ in $E$ as $n \to \infty$ for every $x \in E$. Since $F$ is locally Lipschitz and the $\hat{\Pi}_n$ are uniformly bounded, both as operators from $E$ to $E$ and from $E^*$ to $E^*$, the $F_n$ are locally Lipschitz, where the constant can be chosen uniformly in $n$. Therefore, all of the conditions of Lemma 3.5 are satisfied and we can conclude that $x^n \to x$ in $\mathcal{C}([0,T], E)$ as $n \to \infty$ almost surely.

Using dominated convergence, we see that $P_t^n \varphi(x_n) \to P_t \varphi(x)$ for every $\varphi \in \mathcal{C}_b(E)$ and every $t > 0$, where $(P_t^n)$ are the semigroups from Theorem 3.2 and $(P_t)_{t>0}$ is the semigroup generated by the solutions of (2.1). We can now apply Theorem 3.2 to conclude that $(P_t)_{t>0}$ is $\mu$-symmetric. $\square$

**4. Ergodic properties of the equations.** In this section, we show that the measure $\mu$ from Theorems 3.4 and 3.6 is actually the only invariant measure for both (2.1) and (2.14). This result is essential to justify the use of ergodic averages of solutions to (2.1) or (2.14) in order to sample from $\mu$. We also show that a weak law of large numbers holds for every (and not just almost every) initial condition. Theorems 4.10 and 4.11 summarize the main results.

These results are similar to existing results for (2.1), although our framework includes nonlinear boundary conditions and distribution-valued forcing in the equation. Furthermore, our analysis seems to be completely new for (2.14). The problem is that (2.14) does not have any smoothing property. In particular, it lacks the strong Feller property which is an essential tool in most proofs of uniqueness of invariant measures for SPDEs. We show, however, that it enjoys the recently introduced asymptotic strong Feller property [12], which can, in many cases, be used as a substitute for the strong Feller property, as far as properties of the invariant measures are concerned.

Recall that a Markov semigroup $\mathcal{P}_t$ over a Banach space is called *strong Feller* if it maps bounded measurable functions into bounded continuous functions. It can be shown by a standard density argument that if Assumption 1 holds for $\mathcal{P}_t$, then it also has the strong Feller property. We will not give the precise definition of the asymptotic strong Feller property in the present article since this would require some preliminaries that are not going to be used in the sequel. All we will use is the fact that, in a similar way, if a Markov semigroup $\mathcal{P}_t$ satisfies Assumption 2, then it is also asymptotically strong Feller.

4.1. *Variations of the strong Feller property.* Given a Markov process on a separable Banach space $E$, we call $P_t$ the *associated semigroup* acting on bounded Borel measurable functions $\varphi : E \to \mathbb{R}$. Let us denote by $\mathcal{C}_b^1(E)$ the space of bounded functions from $E$ to $\mathbb{R}$ with bounded Fréchet derivative. For the moment, let us consider processes that satisfy the following property.

ASSUMPTION 1. The Markov semigroup $P_t$ maps $\mathcal{C}_b^1(E)$ into itself. Furthermore, there exists a time $t$ and a locally bounded function $C : E \to \mathbb{R}_+$ such that the bound

$$\|DP_t\varphi(x)\| \leq C(x)\|\varphi\|_\infty \tag{4.1}$$

holds for every $\varphi : E \to \mathbb{R}$ in $\mathcal{C}_b^1(E)$ and every $x \in E$.

It is convenient to introduce

(4.2) $$\mathcal{B}(x) = \{y \in E | \|y - x\|_E \leq 1\}, \qquad \bar{C}(x) = \sup_{y \in \mathcal{B}(x)} C(y).$$

Note that a density argument given in [6] shows that if (4.1) holds for Fréchet differentiable functions, then $P_t \varphi$ is locally Lipschitz continuous with local Lipschitz constant $C(x)\|\varphi\|_\infty$ for every bounded measurable function $\varphi$. In particular, this shows that

(4.3) $$\|P_t(x, \cdot) - P_t(y, \cdot)\|_{\mathrm{TV}} \leq \tfrac{1}{2}\bar{C}(x)\|x - y\|_E$$

for every $x, y \in E$ with $\|x - y\|_E \leq 1$ (with the convention that the total variation distance between mutually singular measures is 1). Recall that the *support* of a measure is the smallest *closed* set with full measure. We also follow the terminology in [6, 23] by calling an invariant measure for a Markov semigroup *ergodic* if the law of the corresponding stationary process is ergodic for the time shifts. The following result follows immediately.

LEMMA 4.1. *Let $P_t$ be a Markov semigroup on a separable Banach space $E$ that satisfies (4.3) and let $\mu$ and $\nu$ be two ergodic invariant measures for $P_t$. If $\mu \neq \nu$, then we have $\|x - y\| \geq \min\{1, 2/\bar{C}(x)\}$ for any two points $(x, y) \in \operatorname{supp} \mu \times \operatorname{supp} \nu$.*

PROOF. Assume (for the purposes of obtaining a contradiction) that there exists a point $(x, y) \in \operatorname{supp} \mu \times \operatorname{supp} \nu$ with $\|x - y\| < 2/\bar{C}(x)$ and $\|x - y\| < 1$. Let $\delta < 1 - \|x - y\|$ be determined later and let $\mathcal{B}_\delta(x)$ denote the ball of radius $\delta$ centered at in $x$. With these definitions, it is easy to check from (4.3) and the triangle inequality that we have

$$\|P_t(x', \cdot) - P_t(y', \cdot)\|_{\mathrm{TV}} \leq \tfrac{1}{2}(2\delta + \|x - y\|)\bar{C}(x)$$

for every $x' \in \mathcal{B}_\delta(x)$ and $y' \in \mathcal{B}_\delta(y)$. Since we assumed that $\|x - y\|\bar{C}(x)/2 < 1$, it is possible, by taking $\delta$ sufficiently small, to find a strictly positive $\alpha > 0$ such that

$$\|P_t(x', \cdot) - P_t(y', \cdot)\|_{\mathrm{TV}} \leq 1 - \alpha.$$

The invariance of $\mu$ and $\nu$ under $P_t$ implies that

$$\|\mu - \nu\|_{\mathrm{TV}} \leq \int_{E^2} \|P_t(\tilde{x}, \cdot) - P_t(\tilde{y}, \cdot)\|_{\mathrm{TV}} \mu(d\tilde{x})\nu(d\tilde{y}) \leq 1 - \alpha\mu(\mathcal{B}_\delta(x))\nu(\mathcal{B}_\delta(y)).$$

Since the definition of the support of a measure implies that both $\mu(\mathcal{B}_\delta(x))$ and $\nu(\mathcal{B}_\delta(y))$ are nonzero, this contradicts the fact that $\mu$ and $\nu$ are distinct and ergodic, therefore mutually singular. $\square$

In our case, it turns out that we are unfortunately not able to prove that (4.1) holds for the equations under consideration. However, it follows immediately from the proof of Lemma 4.1 that we have the following, very similar, result.

COROLLARY 4.2. *Let $P_t$ be a Markov semigroup on a separable Banach space $E$ such that there exists a continuous increasing function $f : \mathbb{R}_+ \to \mathbb{R}_+$ with $f(0) = 0$, $f(1) = 1$ and*

$$\|P_t(x, \cdot) - P_t(y, \cdot)\|_{\mathrm{TV}} \leq \bar{C}(x) f(\|x - y\|) \tag{4.4}$$

*for every $x, y \in E$ with $\|x - y\| \leq 1$. Let $\mu$ and $\nu$ be two ergodic invariant measures for $P_t$. If $\mu \neq \nu$, then we have $f(\|x - y\|) \geq \min\{1, 1/\bar{C}(x)\}$ for any two points $(x, y) \in \operatorname{supp} \mu \times \operatorname{supp} \nu$.*

In Theorem 4.7 below, we will see that the semigroups generated by the nonpreconditioned equations considered in the present article satisfy the smoothing property (4.4). However, even the slightly weaker strong Feller property can be shown to fail for the semigroups generated by the preconditioned equations. They, however, satisfy the following, somewhat weaker, condition.

ASSUMPTION 2. The Markov semigroup $P_t$ maps $\mathcal{C}_b^1(E)$ into itself. Furthermore, there exists a decreasing function $f : \mathbb{R}_+ \to \mathbb{R}_+$ converging to 0 at infinity and a locally bounded function $C : E \to \mathbb{R}_+$ such that the bound

$$\|DP_t \varphi(x)\| \leq C(x)(\|\varphi\|_\infty + f(t)\|D\varphi\|_\infty) \tag{4.5}$$

holds for every $\varphi : E \to \mathbb{R}$ in $\mathcal{C}_b^1(E)$ and every $x \in \mathcal{H}$.

A modification of the argument of Lemma 4.1 yields the following.

LEMMA 4.3. *Let $P_t$ be a Markov semigroup on a separable Banach space $E$ that satisfies Assumption 2 and let $\mu$ and $\nu$ be two ergodic invariant measures for $P_t$. If $\mu \neq \nu$, then we have $\|x - y\| \geq \min\{1, 2/\bar{C}(x)\}$ for any two points $(x, y) \in \operatorname{supp} \mu \times \operatorname{supp} \nu$, where $\bar{C}$ is given in (4.2).*

PROOF. Given a distance $d$ on $E$, recall that the corresponding Wasserstein distance on the space of probability measures on $E$ is given by

$$\|\pi_1 - \pi_2\|_d = \inf_{\pi \in \mathcal{C}(\pi_1, \pi_2)} \int_{E^2} d(x, y) \pi(dx, dy), \tag{4.6}$$

where $\mathcal{C}(\pi_1, \pi_2)$ denotes the set of probability measures on $E^2$ with marginals $\pi_1$ and $\pi_2$.

Given the two invariant measures $\mu$ and $\nu$, we also recall the useful inequality

$$\|\mu - \nu\|_d \leq 1 - \min\{\mu(A), \nu(A)\}\Big(1 - \max_{y,z \in A} \|P_t(z, \cdot) - P_t(y, \cdot)\|_d\Big), \tag{4.7}$$

valid for every $t \geq 0$ and every measurable set $A$ (see, e.g., [12] for a proof).

For $\varepsilon > 0$, we define on $\mathcal{H}$ the distance $d_\varepsilon(x,y) = 1 \wedge \varepsilon^{-1}\|x-y\|$ and denote by $\|\cdot\|_\varepsilon$ the corresponding seminorm on measures given by (4.6). It can be checked from the definitions that, in a way similar to the proof of [12], Proposition 3.12, (4.5) implies that the bound

$$\|P_t(x,\cdot) - P_t(y,\cdot)\|_\varepsilon \leq \frac{1}{2}\|x-y\|\bar{C}(x)\left(1 + \frac{2f(t)}{\varepsilon}\right)$$

holds for every $(x,y) \in E^2$ with $\|x-y\| \leq 1$. Hence, the same argument as in the proof of Lemma 4.1 yields $\alpha > 0$, so, for $\delta$ sufficiently small, we have the bound

$$\|P_t(x,\cdot) - P_t(y,\cdot)\|_\varepsilon \leq (1-\alpha)\left(1 + \frac{2f(t)}{\varepsilon}\right)$$

for every $x' \in \mathcal{B}_\delta(x)$ and $y' \in \mathcal{B}_\delta(y)$. Note that $\delta$ can be chosen independently of $\varepsilon$. Choosing $t$ as a function of $\varepsilon$ sufficiently large so that $f(t) < \alpha\varepsilon/2$, say, it follows from (4.7) that

$$\|\mu_1 - \mu_2\|_\varepsilon \leq 1 - \alpha^2 \min\{\mu_1(\mathcal{B}_\delta(x)), \mu_2(\mathcal{B}_\delta(x))\}$$

for every $\varepsilon > 0$. Since $\lim_{\varepsilon \to 0}\|\mu_1 - \mu_2\|_\varepsilon = \|\mu_1 - \mu_2\|_{\mathrm{TV}}$ (see [12]), the claim follows in the same way as in the proof of Lemma 4.1. $\square$

4.2. *Conditions for* (4.4) *to hold.* In this subsection, we show that equation (2.1), arising from the nonpreconditioned case, satisfies the bound (4.4). Our main result is the following theorem.

The proof of the results is closely related to standard arguments that can be found, for example, in [4, 6, 17]. However, the situation in these works is different from ours, mainly because we only have local bounds on the derivative of the flow with respect to the initial condition. This forces us to use an approximation argument which, in turn, only yields a bound of type (4.4) rather than the bound (4.1) obtained in the previously mentioned works. The present proof unfortunately requires (4.8) as an additional assumption on the nonlinearity $F$, even though we believe that this is somewhat artificial.

THEOREM 4.4. *Suppose that assumptions* (A1)–(A4) *hold. Assume, furthermore, that for every $R > 0$, there exists a Fréchet differentiable function $F_R : E \to E^*$ such that*

(4.8)
$$\begin{aligned} F_R(x) &= F(x) & &\text{for } \|x\|_E \leq R, \\ F_R(x) &= 0 & &\text{for } \|x\|_E \geq 2R \end{aligned}$$

*and such that there exist constants $C$ and $N$ such that*

$$\|F_R(x)\| + \|DF_R(x)\| \leq C(1+R)^N$$

*for every* $x \in E$. *There then exist exponents* $\tilde{N} > 0$ *and* $\alpha > 0$ *such that the solutions to the SPDE* (2.1) *satisfy* (4.4) *with* $f(r) = r^\alpha$ *and* $\bar{C}(x) = (1 + \|x\|_E)^{\tilde{N}}$.

PROOF. Fix $x_0 \in E$ and define $R = 2\|x_0\|_E$. Denote by $\Phi_t^R : E \to E$ the flow induced by the solutions to the truncated equation

$$(4.9) \qquad dx = \mathcal{L}x \, dt + F_R(x) \, dt + \sqrt{2} \, dw(t).$$

Further, denote by $z$ the solution to the linearized equation defined in (2.4). It follows immediately from Picard iterations that $\Phi_t^R$ is Fréchet differentiable and that there exists a constant $C$ such that

$$(4.10) \quad \begin{aligned} \|\Phi_t^R(x)\|_E &\leq \|x\|_E + \|z(t)\|_E + Ct^{1-2\alpha}(1+R)^N, \\ \|D\Phi_t^R(x)\|_{E \to E} &\leq 2 \end{aligned}$$

for every $t$ with

$$(4.11) \qquad t^{1-2\alpha} \leq \frac{1}{C(1+R)^N}.$$

Note that the bounds in (4.10) are almost sure bounds and that (4.11) is a *deterministic* condition on the time interval we are allowed to consider.

Now, denote by $P_t^R$ the Markov semigroup generated by (4.9). For an arbitrary function $\varphi \in \mathcal{C}_b^1(E)$ and an arbitrary vector $\xi \in E$, the Bismut–Elworthy–Li formula [6, 9] yields

$$|DP_t^R \varphi(x) \xi| = \frac{1}{t} \mathbb{E} \left( \varphi(\Phi_t^R(x)) \int_0^t \langle D\Phi_s^R(x)\xi, dw(s) \rangle \right)$$

$$\leq \frac{1}{t} \|\varphi\|_\infty \left( \mathbb{E} \int_0^t \|D\Phi_s^R(x)\xi\|_{\mathcal{H}}^2 \, ds \right)^{1/2}.$$

Combining this with (4.10) shows that there exists a constant $C$ such that

$$(4.12) \qquad \|P_t^R(x, \cdot) - P_t^R(y, \cdot)\|_{\mathrm{TV}} \leq \frac{C}{\sqrt{t}} \|x - y\|_E,$$

provided that $t$ is sufficiently small that (4.11) holds. The bound (4.10) shows that there exists $\theta > 0$ such that

$$(4.13) \qquad \mathbb{P}\left( \sup_{s \in [0,t]} \|x(s)\|_E \geq R \right) \leq \frac{Ct^\theta}{R}$$

for every $t$ such that (4.11) holds and every $x_0$ such that $\|x_0\|_E \leq R/2$.

Furthermore, it is clear that the solution to (4.9) agrees with the solution to (2.1), provided it stays inside a ball of radius $R$, so (4.13) implies that, under the same conditions, we have

$$(4.14) \qquad \|P_t(x, \cdot) - P_t^R(x, \cdot)\|_{\mathrm{TV}} \leq \frac{Ct^\theta}{R}.$$

Combining (4.14) and (4.12) yields

$$\|P_t(x,\cdot) - P_t(y,\cdot)\|_{\mathrm{TV}} \leq \frac{C}{\sqrt{t}}\|x-y\|_E + \frac{Ct^\theta}{R} \tag{4.15}$$

for all pairs $(x,y) \in E \times E$ such that $\sup\{\|x\|_E, \|y\|_E\} \leq R/2$ and all times $t$ satisfying (4.11). Since we have

$$\|P_t(x,\cdot) - P_t(y,\cdot)\|_{\mathrm{TV}} \leq \|P_s(x,\cdot) - P_s(y,\cdot)\|_{\mathrm{TV}}$$

for $s \leq t$, (4.15) actually implies that

$$\|P_t(x,\cdot) - P_t(y,\cdot)\|_{\mathrm{TV}} \leq \inf_{s \leq t}\left(\frac{C}{\sqrt{s}}\|x-y\|_E + \frac{Cs^\theta}{R}\right),$$

which immediately yields that a bound of the type (4.4) holds, with $\bar{C}(x)$ growing polynomially in $\|x\|_E$. □

COROLLARY 4.5. *Let $U : E \to \mathbb{R}$ be bounded from above and Fréchet differentiable. Assume that $\mathcal{L}$ and $F = U'$ satisfy the assumptions of Theorem 4.4. The SDE* (2.1) *then has a unique stationary distribution, which is given by* (3.7).

PROOF. Denote by $\mathcal{E}$ the set of all ergodic invariant measures for (2.1). It follows from Theorem 3.4 that $\mu$, as given by (3.7), is an invariant measure for (2.1), so $\mathcal{E}$ is not empty. Also, note that the support of $\mu$ is equal to $E$ since the embedding $\mathcal{H}^{1/2} \subset E$ is dense by (A2). Assume, now, that $\mathcal{E}$ contains at least two elements, $\nu_1$ and $\nu_2$. In this case, it follows from Theorem 4.4 that there exists an open set $A \subset E$ such that $A \cap \operatorname{supp}\nu = \varnothing$ for every $\nu \in \mathcal{E}$. Since every invariant measure is a convex combination of ergodic invariant measures ([23], Theorem 6.6), this implies that $\mu(A) = 0$, which contradicts to the fact that $\operatorname{supp}\mu = E$. □

REMARK 4.6. Since we obtain the strong Feller property for (2.1), as well as the existence of a Lyapunov function [see equation (2.16)], we can apply the machinery exposed in [18] in order to obtain the exponential convergence (in a weighted total variation norm) of transition probabilities to the unique invariant measure. The only additional ingredient that is required is the fact that the level sets of the Lyapunov function are "small." This can be checked by a standard controllability argument.

4.3. *Conditions for* (4.5) *to hold.* In this subsection, we show that the equations arising from the non-preconditioned case satisfy a bound of the type (4.5).

THEOREM 4.7. *Let $\mathcal{L}$, $F$ and $\mathcal{G}$ satisfy* (A1)–(A3) *and* (A5)–(A7). *The Markov semigroup on $\mathcal{H}$ generated by the solutions of* (2.14) *then satisfies the bound* (4.5), *with $C(x) \leq C(1 + \|x\|_E)^N$ for some constants $C$ and $N$. In particular, it is asymptotically strong Feller.*

REMARK 4.8. Note that it is not generally true that these assumptions imply that the process is strong Feller. A counterexample is given by the case where $\tilde{\mathcal{L}}$ is minus the identity, $F = 0$ and $\mathcal{G} : \mathcal{H} \to \mathcal{H}$ is any positive definite trace class operator. This counterexample comes very close to the situation studied in this paper, so the strong Feller property is clearly not an appropriate concept here.

PROOF OF THEOREM 4.7. It follows from standard arguments that the evolution map $\Phi_{s,t} : E \times \Omega \to E$ is Fréchet differentiable. In the sequel, we denote its Fréchet derivative by $J_{s,t}$.

The family of (random) linear operators $J_{s,t} : E \to E$ is given in the following way. For every $\xi \in E$, $J_{s,t}\xi$ solves the equation

$$\partial_t J_{s,t}\xi = \tilde{\mathcal{L}} J_{s,t}\xi + D\tilde{F}(x(t)) J_{s,t}\xi, \qquad J_{s,s}\xi = \xi.$$

We also define a family of (random) linear operators $A_t : \mathrm{L}^2([0,t], \mathcal{H}) \to E$ by

$$A_t v = \int_0^t J_{s,t} \mathcal{G}^{1/2} v(s) \, ds.$$

This is well defined since $\mathcal{G}^{1/2}$ maps $\mathcal{H}$ into $E$ by Lemma 2.8. Recall that $A_t v$ is the Malliavin derivative of the flow at time $t$ in the direction of the Cameron–Martin vector $v$. We will also denote this by $A_t v = \mathcal{D}^v \Phi_{0,t}$.

Given a perturbation $\xi$ in the initial condition for $x$, the idea is to find a perturbation $v$ in the direction of the Cameron–Martin space of the noise such that these perturbations "cancel" each other for large times $t$. Given a square-integrable $\mathcal{H}$-valued process $v$, we therefore introduce the notation

$$\rho(t) = J_{0,t}\xi - A_t v_{[0,t]},$$

where $v_J$ denotes the restriction of $v$ to the interval $J$. Note that $\rho(t)$ is the solution to the differential equation

$$(4.16) \qquad \partial_t \rho(t) = \tilde{\mathcal{L}}\rho(t) + D\tilde{F}(x(t))\rho(t) - \mathcal{G}^{1/2} v(t), \qquad \rho(0) = \xi \in E.$$

The reason for introducing this process $\rho$ is clear from the approximate integration by parts formula (see [12] for more details), which holds for every bounded function $\varphi : E \to \mathbb{R}$ with bounded Fréchet derivative:

$$\langle DP_t \varphi(x), \xi \rangle = \mathbb{E}(\langle D(\varphi(x_t)), \xi \rangle) = \mathbb{E}((D\varphi)(x_t) J_{0,t}\xi)$$

$$\begin{aligned}
&= \mathbb{E}((D\varphi)(x_t)A_t v_{[0,t]}) + \mathbb{E}((D\varphi)(x_t)\rho_t) \\
&= \mathbb{E}(\mathcal{D}^{v_{[0,t]}}\varphi(x_t)) + \mathbb{E}((D\varphi)(x_t)\rho_t) \\
&= \mathbb{E}\bigg(\varphi(x_t)\int_0^t \langle v(s), dw(s)\rangle\bigg) + \mathbb{E}((D\varphi)(x_t)\rho_t) \\
&\leq \|\varphi\|_\infty \sqrt{\mathbb{E}\int_0^t \|v(s)\|^2\, ds} + \|D\varphi\|_\infty \mathbb{E}\|\rho_t\|_E.
\end{aligned} \quad (4.17)$$

In this formula, $w$ denotes a cylindrical Wiener process on $\mathcal{H}$, so $\tilde{w} = \mathcal{G}^{1/2}w$. This formula is valid for every adapted square integrable $\mathcal{H}$-valued process $v$.

It remains to choose an adapted process $v$ such that $\rho(t) \to 0$. For

$$v(t) = \mathcal{G}^{-1/2}(D\tilde{F}(x(t)) + \mathcal{K})e^{-t}\xi,$$

it is easy to check that equation (4.16) reduces to $\partial_t \rho = -\rho$, so $\|\rho(t)\|_E = e^{-t}$. Furthermore, Theorem 2.10, together with assumptions (A3) and (A7), ensures that $\mathbb{E}\|v(t)\|_\mathcal{H}^2 \leq C(1+\|x\|_E)^N e^{-wt}$ for some constants $C$, $N$ and $w$, so (4.17) immediately implies (4.5). □

COROLLARY 4.9. *Let $U: E \to \mathbb{R}$ be bounded from above, Fréchet differentiable and such that (A1)–(A3) and (A5)–(A7) hold for $F = U'$. The SDE (2.14) then has a unique stationary distribution, given by (3.11).*

PROOF. The proof follows exactly the same pattern as the proof of 4.5, but we replace references to Theorem 4.4 by references to Theorem 4.7 and 4.3. □

4.4. *Law of large numbers.* In this section, we use the results of the previous section in order to show that the solutions to our equations satisfy a law of large numbers. We first introduce the following result.

THEOREM 4.10. *Assume that (A1)–(A4) and (4.8) hold and let $\mu$ be an ergodic invariant probability measure for (2.1). We then have*

$$\lim_{T\to\infty} \frac{1}{T}\int_0^T \varphi(x(t))\, dt = \int_E \varphi(x)\mu(dx) \qquad \text{almost surely} \quad (4.18)$$

*for every initial condition $x_0$ in the support of $\mu$ and for every bounded measurable function $\varphi: E \to \mathbb{R}$.*

PROOF. Denote by $\mathcal{A} \subset E$ the set of initial conditions for which (4.18) holds and by $S$ the support of $\mu$. We know from Birkhoff's ergodic theorem that $\mu(\mathcal{A}) = 1$ and therefore that $\mathcal{A}$ is dense in $S$. Now, let $x_0 \in S$ and $\varepsilon > 0$ be arbitrary and choose a sequence of points $x_0^n$ in $\mathcal{A}$ converging to $x_0$.

Fix an arbitrary time $t_0 > 0$. Since, by Theorem 4.4, $P_t(x, \cdot)$ is continuous in $x$ (in the topology of total variation), there exists $n$ such that

$$\|P_{t_0}(x_0^n, \cdot) - P_{t_0}(x_0, \cdot)\|_{\mathrm{TV}} < \varepsilon. \tag{4.19}$$

Let $x^n(\cdot)$ denote the trajectories starting from $x_0^n$ and $x(t)$ denote the trajectories starting from $x_0$. By the Markov property, the bound (4.19) implies that there exists a coupling between the laws of $x^n(\cdot)$ and $x(\cdot)$ such that, with probability larger than $1 - \varepsilon$, we have $x^n(t) = x(t)$ for every $t \geq t_0$. This immediately shows that

$$\lim_{T \to \infty} \frac{1}{T} \int_0^T \varphi(x(t))\, dt = \int_E \varphi(x) \mu(dx)$$

on a set of measure larger than $1 - \varepsilon$. Since $\varepsilon$ was arbitrary, the desired result follows. □

In the preconditioned case, we have the following, somewhat weaker, form of the law of large numbers.

THEOREM 4.11. *Assume that* (A1)–(A3) *and* (A5)–(A7) *hold and let $\mu$ be an ergodic invariant probability measure for* (2.14). *We then have*

$$\lim_{T \to \infty} \frac{1}{T} \int_0^T \varphi(x(t))\, dt = \int_E \varphi(x) \mu(dx) \qquad \text{in probability} \tag{4.20}$$

*for every initial condition $x(0)$ in the support of $\mu$ and for every bounded function $\varphi : E \to \mathbb{R}$ with bounded Fréchet derivative.*

PROOF. As before, denote by $\mathcal{A} \in E$ the set of initial conditions for which (4.20) holds and by $S$ the support of $\mu$. Since convergence in probability is weaker than almost sure convergence, we know from Birkhoff's ergodic theorem that $\mu(\mathcal{A}) = 1$ and therefore that $\mathcal{A}$ is dense in $S$.

Define

$$\mathcal{E}_\varphi^T(x) = \frac{1}{T} \int_0^T \varphi(x(t))\, dt \qquad \text{with } x(0) = x.$$

The idea is to use the following chain of equalities, valid for every pair of bounded functions $\varphi : E \to \mathbb{R}$ and $\psi : \mathbb{R} \to \mathbb{R}$ with bounded Fréchet derivatives. The symbol $D$ denotes the Fréchet derivative of a given function and the symbol $\mathcal{D}$ denotes its Malliavin derivative. We have

$$D\mathbb{E}\psi(\mathcal{E}_\varphi^T(x))\xi = \mathbb{E}\left((D\psi)(\mathcal{E}_\varphi^T)\frac{1}{T}\int_0^T (D\varphi)(x(t)) J_{0,t}\xi\, dt\right)$$

$$= \mathbb{E}(\mathcal{D}^v \psi(\mathcal{E}_\varphi^T)) + \mathbb{E}\left((D\psi)(\mathcal{E}_\varphi^T)\frac{1}{T}\int_0^T (D\varphi)(x(t))\rho(t)\, dt\right)$$

$$\leq \mathbb{E}\left(\psi(\mathcal{E}_\varphi^T)\int_0^T v(t)\,dt\right) + \frac{\|D\psi\|_\infty \|D\varphi\|_\infty}{T}\mathbb{E}\int_0^T |\rho(t)|\,dt$$

$$\leq C\left(\|\psi\|_\infty + \frac{\|D\psi\|_\infty \|D\varphi\|_\infty}{T}\right)\|\xi\|.$$

Now, denote by $\mu_\varphi^T(x)$ the law of $\mathcal{E}_\varphi^T(x)$. The above chain of inequalities shows that

$$\|\mu_\varphi^T(x) - \mu_\varphi^T(y)\|_W \leq C\left(1 + \frac{\|D\varphi\|_\infty}{T}\right)\|x - y\|$$

for some constant $C$, where $\|\cdot\|_W$ denotes the Wasserstein distance between two probability measures with respect to the distance function $1 \wedge \|x - y\|$. Since the Wasserstein distance metrizes the weak convergence topology and weak convergence to a delta measure is the same as convergence in probability to the point at which the measure is located, this implies that (4.20) holds for *every* initial condition $x$ in $S$. $\square$

REMARK 4.12. It is possible to extend the above argument to a larger class of continuous test functions $\varphi$ by introducing a time-dependent smoothing (and possibly cu-toff).

REMARK 4.13. If we wish to obtain a statement which is valid for *every* initial condition, it is, in general, impossible to drop the continuity assumption on $\varphi$. Consider, for example, the trivial dynamic $\dot{x} = -x$ on $\mathbb{R}$ with invariant measure $\delta_0$. It is obvious that if we take $x_0 = 1$, $\varphi(0) = 1$ and $\varphi(x) = 0$ for $x \neq 0$, then the left-hand side of (4.20) is 0, whereas the right-hand side is 1.

**5. Conditioned SDEs.** In this section, we outline how the preceding material can be used to construct SPDEs which sample from the distribution of conditioned SDEs. The program outlined here will be carried out in the subsequent sections in three specific contexts.

We start the section by explaining the common structure of the arguments used in each of the following three sections; we also outline the common technical tools required. We then make some remarks concerning the conversion between Hilbert-space-valued SDEs and SPDEs and, in particular, discuss how the framework developed in preceding sections enables us to handle the nonlinear boundary conditions which arise.

Consider the following $\mathbb{R}^d$-valued SDEs, both driven by a $d$-dimensional Brownian motion, with invertible covariance matrix $BB^*$:

(5.1) $$dX = AX\,du + f(X)\,du + B\,dW, \qquad X(0) = x^-$$

and

(5.2) $$dZ = AZ\,du + B\,dW, \qquad Z(0) = x^-.$$

Our aim is to construct an SPDE which has the distribution of $X$, possibly conditioned by observations, as its stationary distribution. The construction consists of the following steps. We symbolically denote the condition on $X$ and $Z$ by C here and we set $m(u) = \mathbb{E}(Z(u)|\mathtt{C})$.

1. Use the Girsanov formula (Lemma 5.2 below) to find the density of the distribution $\mathcal{L}(X)$ w.r.t. $\mathcal{L}(Z)$.

2. Use results about conditional distributions (Lemma 5.3 below) to derive the density of the conditional distribution $\mathcal{L}(X|\mathtt{C})$ w.r.t. $\mathcal{L}(Z|\mathtt{C})$. Using substitution, this gives the density of the shifted distribution $\mathcal{L}(X|\mathtt{C}) - m$ w.r.t. the centered measure $\mathcal{L}(Z|\mathtt{C}) - m$.

3. Use the results of the companion paper [14] to obtain an $L^2$-valued SDE which has the centered Gaussian measure $\mathcal{L}(Z|\mathtt{C}) - m$ as its stationary distribution. This also gives a representation of $m$ as the solution of a boundary value problem.

4. Use the results of Sections 2 and 3 and the density from step 2 to derive an $C([0,1],\mathbb{R}^d)$-valued SDE with stationary distribution $\mathcal{L}(X|\mathtt{C}) - m$. Use the results of Section 4 to show ergodicity of the resulting SDE.

5. Write the $L^2$-valued SDE as an SPDE, reversing the centering from step 2 in the process.

Combining all of these steps leads to an SPDE which samples from the conditional distribution $\mathcal{L}(X|\mathtt{C})$ in its stationary measure. In the remaining part of this section, we will elaborate on the parts of the outlined program which are common to all three of our applications.

We will assume throughout the rest of this article that the drift $f$ for $X$ is of the form $f = -BB^*\nabla V$, where the potential $V$ satisfies the following polynomial growth condition.

(M) The potential $V: \mathbb{R}^d \to \mathbb{R}$ is a $\mathcal{C}^4$-function which can be written as

$$V(x) = M(x, \ldots, x) + \tilde{V}(x),$$

where $M: (\mathbb{R}^d)^{2p} \to \mathbb{R}$ is $2p$-linear with

$$M(x, \ldots, x) \geq c|x|^{2p} \qquad \forall x \in \mathbb{R}^d$$

for some $p \in \mathbb{N}$ and $c > 0$, and $\tilde{V}: \mathbb{R}^d \to \mathbb{R}$ satisfies

$$\frac{|D^k \tilde{V}(x)|}{1 + |x|^{2p-k}} \to 0 \qquad \text{as } |x| \to \infty$$

for every $k$-fold partial derivative operator $D^k$ with $k = 0, \ldots, 4$.

Under condition (M), the potential $V$ is bounded from below and grows like $|x|^{2p}$ as $|x| \to \infty$. From [16], Section 2.3, Theorem 3.6, we know that under this condition on $f$, the SDE (5.1) has a unique nonexploding solution.

Later, when checking assumption (A4) and the boundedness of $U$ in Theorem 3.4, we have to estimate terms which involve both the nonlinearity $f$ and the linear part $A$ of the drift. If condition (M) is satisfied for $p > 1$, we will get the estimates from the superlinear growth of $f$. For $p = 1$, we use the following, additional assumption on $A$.

(Q) For $p = 1$, the matrices $A, B$ from (5.1) satisfy $QA + A^*Q - QBB^*Q < 0$ (as a symmetric matrix), where $Q \in \mathbb{R}^{d \times d}$ is the symmetric matrix defined by the relation $M(x,x) = \frac{1}{2}\langle x, Qx \rangle$ for all $x \in \mathbb{R}^d$.

NOTATION 5.1. Introduce the inner product and related norm

$$\langle a, b \rangle_B = a^*(BB^*)^{-1}b, \qquad |a|_B^2 = \langle a, a \rangle_B,$$

defined for any invertible $B$.

The densities in step 1 above will be calculated from the Girsanov formula. As an abbreviation, let

(5.3) $$\Phi = \tfrac{1}{2}(|f|_B^2 + \operatorname{div} f).$$

When expressed in terms of $V$, this becomes

$$\Phi = \tfrac{1}{2}(|B^*\nabla V|^2 - (BB^*):D^2V),$$

where ":" denotes the Frobenius inner product and $D^2V$ denotes the Hessian of $V$.

LEMMA 5.2. *Assume that* (5.1) *has a solution without explosions on the interval* $[0,1]$. *Let* $\mathbf{Q}$ *(resp.* $\mathbf{P}$*) be the distribution on path space* $\mathcal{C}([0,1], \mathbb{R}^d)$ *of the solution of* (5.2) *[resp.* (5.1)*]. Then,*

$$d\mathbf{P}(Z) = \frac{1}{\varphi(Z)} d\mathbf{Q}(Z),$$

*where*

$$\ln \varphi(Z) = -\int_0^1 \langle f(Z(u)), \circ dZ(u) \rangle_B + \int_0^1 \left( \Phi(Z(u)) + \langle f(Z(u)), AZ(u) \rangle_B \right) du.$$

PROOF. Since $X$ (by assumption) and $Z$ (since it solves a linear SDE) have no explosions, we can apply Girsanov's theorem ([8], Theorem 11A) which yields

$$\ln \varphi(Z) = -\int_0^1 \langle B^{-1}f(Z(u)), dW(u) \rangle - \int_0^1 \tfrac{1}{2}|f(Z(u))|_B^2 \, du.$$

But,

$$\int_0^1 \langle B^{-1} f(Z(u)), dW(u) \rangle$$
$$= \int_0^1 \langle f(Z(u)), dZ(u) \rangle_B - \langle f(Z(u)), AZ(u) \rangle_B \, du - \int_0^1 |f(Z(u))|_B^2 \, du.$$

Converting the first integral on the right-hand side to Stratonovich form gives the desired result. $\square$

Writing the Radon–Nikodym derivative in terms of a Stratonovich integral in the lemma above is helpful when studying its form in the case of gradient vector fields; the stochastic integral then reduces to boundary contributions.

We will handle the conditioning in step 2 of the program outlined above with the help of Lemma 5.3 below. We will use it in two ways: to condition on paths $X$ which end at $X(1) = x^+$ and, for the filtering/smoothing problem where $X$ will be replaced by a pair $(X, Y)$, to condition the signal $(X(u))_{u \in [0,1]}$ on the observation $(Y(u))_{u \in [0,1]}$. Since the proof of the lemma is elementary, we omit it here (see Section 10.2 of [7] for reference).

LEMMA 5.3. *Let $\mathbf{P}, \mathbf{Q}$ be probability measures on $S \times T$, where $(S, \mathcal{A})$ and $(T, \mathcal{B})$ are measurable spaces, and let $X : S \times T \to S$ and $Y : S \times T \to T$ be the canonical projections. Assume that $\mathbf{P}$ has a density $\varphi$ w.r.t. $\mathbf{Q}$ and that the conditional distribution $\mathbf{Q}_{X|Y=y}$ exists. The conditional distribution $\mathbf{P}_{X|Y=y}$ then exists and is given by*

(5.4) $$\frac{d\mathbf{P}_{X|Y=y}}{d\mathbf{Q}_{X|Y=y}}(x) = \begin{cases} \dfrac{1}{c(y)} \varphi(x, y), & \text{if } c(y) > 0, \\ 1, & \text{otherwise,} \end{cases}$$

*with $c(y) = \int_S \varphi(x, y) \, d\mathbf{Q}_{X|Y=y}(x)$ for all $y \in T$.*

The linear, infinite-dimensional SDEs from [14] which we will use in step 3 are defined on the space $\mathcal{H} = \mathrm{L}^2([0,1], \mathbb{R}^d)$ and the generator of the corresponding semigroup is the self-adjoint operator $\mathcal{L}$ on $\mathcal{H}$ which is the extension of the differential operator

$$L = (\partial_u + A^*)(BB^*)^{-1}(\partial_u - A)$$

with appropriate boundary conditions. When studying the filtering/smoothing problem, the operator $L$ will include additional lower order terms, which we omit here for clarity.

The nonlinear, infinite-dimensional SDEs derived in step 4 are of the form (2.1) or (2.14). They share the operator $\mathcal{L}$ with the linear equations, but have an additional nonlinear drift $F : E \to E^*$, where the space $E$ will be a subspace of $\mathcal{C}([0,1], \mathbb{R}^d)$. The main difficulty in step 4 is to verify that assumptions (A1)–(A4) for the nonpreconditioned case or (A1)–(A3), (A5)–(A7) for the preconditioned case hold under conditions (M) and (Q). The nonlinearity $F$ is of the form

$$(5.5) \qquad (F(\omega))(u) = \varphi(\omega(u)) + h_0(\omega(0))\delta(u) + h_1(\omega(1))\delta(1-u)$$

for all $u \in [0,1]$, where $\varphi$, $h_0$ and $h_1$ are functions from $\mathbb{R}^d$ to $\mathbb{R}^d$. The symbols $\delta(u)$ and $\delta(1-u)$ denote Dirac mass terms at the boundaries. The functions $\varphi$, $h_0$ and $h_1$ are calculated from the potential $V$ and, in our applications, the growth conditions from (A3) will be a direct consequence of condition (M). The following lemma, in conjunction with condition (M), will help us to verify assumption (A4).

LEMMA 5.4. *Let $c, \gamma > 0$ and $h : \mathbb{R}^d \to \mathbb{R}^d$ be continuous with $\langle h(x), x \rangle \leq -\gamma |x|^2$ for every $x \in \mathbb{R}^d$ with $|x| > c$. Then,*

$$\langle \omega^*, h(\omega) \rangle \leq -\gamma \|\omega\|_\infty \qquad \forall \omega^* \in \partial \|\omega\|_\infty$$

*and for all continuous functions $\omega : [0,1] \to \mathbb{R}^d$ such that $\|\omega\|_\infty \geq c$.*

PROOF. Using the characterization (2.8) of $\partial \|\omega\|_\infty$ from the remark after (A4), we get

$$\langle \omega^*, h(\omega) \rangle = \int_0^1 \left\langle h(\omega(u)), \frac{\omega(u)}{|\omega(u)|} \right\rangle |\omega^*|(du)$$

$$\leq -\gamma \int_0^1 |\omega(u)| |\omega^*|(du) = -\gamma \|\omega\|_\infty.$$

This completes the proof. □

REMARK 5.5. The special choice of $E$, $\mathcal{L}$ and $F$ allows us to rewrite the Hilbert-space-valued SDEs as $\mathbb{R}^d$-valued SPDEs in step 5. We obtain SPDEs of the following form:

$$\partial_t x(t,u) = Lx(t,u) + g(u) + \varphi(x(t,u)) + h_0(x(t,0))\delta(u) + h_1(x(t,1))\delta(1-u)$$
$$+ \sqrt{2} \partial_t w(t,u) \qquad \forall (t,u) \in (0,\infty) \times [0,1],$$

$$D_0 x(t,0) = \alpha, \qquad D_1 x(t,1) = \beta \qquad \forall t \in (0,\infty),$$

where $\varphi, h_0, h_1$ are functions from $\mathbb{R}^d$ to $\mathbb{R}^d$, $\partial_t w$ is space-time white noise, $D_i = A_i \partial_u + B_i$ are linear first-order differential operators and $\alpha, \beta \in \mathbb{R}^d$ are constants. The term $g$ is only nonzero for the filtering/smoothing problem and is then an element of $E^*$. Incorporating the jump induced by the Dirac masses into the boundary conditions gives

$$D_0 x(t,0) = \alpha - A_0(BB^*)h_0(x(t,0)),$$
$$D_1 x(t,1) = \beta + A_1(BB^*)h_1(x(t,1)) \qquad \forall t \in (0,\infty).$$

With these boundary conditions, the delta functions are removed from the SPDE above.

We call a process $x: [0,\infty) \times [0,1] \to \mathbb{R}^d$ a *mild solution* of this SPDE if $x - m$ is a mild solution of the $\mathcal{H}$-valued SDE (2.1), where $m$ is a solution of the boundary value problem $-Lm = g$ with $D_0 m(0) = \alpha$ and $D_1 m(1) = \beta$ and $\mathcal{L}$ is the self-adjoint operator $L$ with boundary conditions $D_0 \omega(0) = 0$ and $D_1 \omega(1) = 0$.

REMARK 5.6. When using the preconditioned equation (2.14), we will consider evolution equations of the following form:

$$\partial_t x(t,u) = -x(t,u) + y(t,u) + \sqrt{2}\partial_t \tilde{w}(t,u) \qquad \forall (t,u) \in (0,\infty) \times [0,1],$$
$$-L_0 y(t,u) = L_1 x(t,u) + g(u) + \varphi(x(t,u))$$
$$\qquad + h_0(x(t,0))\delta(u) + h_1(x(t,1))\delta(1-u)$$
$$\qquad\qquad\qquad\qquad\qquad \forall (t,u) \in (0,\infty) \times [0,1],$$

$$D_0 y(t,0) = \alpha, \qquad D_1 y(t,1) = \beta \qquad \forall t \in (0,\infty),$$

where $L = L_0 + L_1$, $L_0$ is a second-order differential operator, $L_1$ is a differential operator of lower order, $\mathcal{G}$ is the inverse of $-L_0$ subject to the same homogeneous boundary conditions as $L$ and $\tilde{w}$ is a $\mathcal{G}$-Wiener process. Incorporating the induced jump into the boundary conditions as above gives

$$D_0 y(t,0) = \alpha - A_0(BB^*)h_0(x(t,0)),$$
$$D_1 y(t,1) = \beta + A_1(BB^*)h_1(x(t,1)) \qquad \forall t \in (0,\infty).$$

With these boundary conditions, the Dirac mass is removed from the evolution equation above.

We call a process $x: [0,\infty) \times [0,1] \to \mathbb{R}^d$ a *strong solution* of this SPDE if $x - m$ is a strong solution of the $\mathcal{H}$-valued SDE (2.14), where $m$ is a solution of the boundary value problem $-L_0 m = g$ with $D_0 m(0) = \alpha$ and $D_1 m(1) = \beta$ and $\mathcal{L}$ is the self-adjoint operator $L = L_0 + L_1$ with boundary conditions $D_0 \omega(0) = 0$ and $D_1 \omega(1) = 0$.

**6. Free path sampling.** In this section, we will follow the program outlined in Section 5 in order to construct SPDEs whose stationary distribution is the distribution of the solution $X$ of the SDE (5.1). The main results are Theorems 6.1 and 6.3. We re-emphasize that it is straightforward to generate independent samples from the desired distribution in this unconditioned case and there would be no reason to use the SPDEs in practice for this problem. However, the analysis highlights a number of issues which arise in the two following sections, in a straightforward way; we therefore include it here.

We write $\mathcal{C}_-([0,1],\mathbb{R}^d)$ for the set of all continuous functions from $[0,1]$ to $\mathbb{R}^d$ with start at $x^-$.

THEOREM 6.1. *Let $A \in \mathbb{R}^{d \times d}$ be a matrix, let $B \in \mathbb{R}^{d \times d}$ be invertible, let $f = -BB^*\nabla V$, assume that conditions* (M) *and* (Q) *are satisfied and let $x^- \in \mathbb{R}^d$. Consider the $\mathbb{R}^d$-valued SPDE*

$$\partial_t x = (\partial_u + A^*)(BB^*)^{-1}(\partial_u - A)x - \nabla \Phi(x) \tag{6.1a}$$
$$- (Df)^*(x)(BB^*)^{-1}Ax - A^*(BB^*)^{-1}f(x) + \sqrt{2}\partial_t w,$$

(6.1b) $\quad x(t,0) = x^-, \qquad \partial_u x(t,1) = Ax(t,1) + f(x(t,1)),$

(6.1c) $\quad x(0,u) = x_0(u),$

*where $\partial_t w$ is space-time white noise and $\Phi$ is given by (5.3).*

(a) *This SPDE has a unique, mild solution for every $x_0 \in \mathcal{C}_-([0,1],\mathbb{R}^d)$ and its stationary distribution coincides with the distribution of the solution of SDE (5.1).*

(b) *For every bounded, measurable function $\varphi : \mathcal{C}_-([0,1],\mathbb{R}^d) \to \mathbb{R}$ and every $x_0 \in \mathcal{C}_-([0,1],\mathbb{R}^d)$, we have*

$$\lim_{T \to \infty} \frac{1}{T} \int_0^T \varphi(x(t,\cdot))\,dt = \mathbb{E}(\varphi(X)) \qquad \textit{almost surely,}$$

*where $X$ is the solution of (5.1).*

PROOF. Let $X$ be a solution of (5.1) and let $Z$ be the solution of the linear SDE (5.2). From Lemma 5.2, we know that the distribution of $X$ has a density $\varphi$ with respect to the distribution of $Z$ which is given by

$$\varphi(\omega) = c \cdot \exp\Bigl(-V(\omega(1)) \tag{6.2}$$
$$+ \int_0^1 \langle \nabla V(\omega(u)), A\omega(u) \rangle \,du - \int_0^1 \Phi(\omega(u))\,du \Bigr)$$

for all $\omega \in \mathcal{C}([0,1],\mathbb{R}^d)$ and some normalization constant $c$. Let $m(u) = \mathbb{E}(Z(u))$ for all $u \in [0,1]$. The density $\psi$ of the distribution $\mu = \mathcal{L}(X - m)$

w.r.t. the centered distribution $\nu = \mathcal{L}(Z - m)$ is then given by $\psi(\omega - m) = \varphi(\omega)$ for all $\omega \in \mathcal{C}([0,1], \mathbb{R}^d)$.

Consider the Hilbert space $\mathcal{H} = L^2([0,1], \mathbb{R}^d)$ and the Banach space $E = \{\omega \in \mathcal{C}([0,1], \mathbb{R}^d) | \omega(0) = 0\} \subseteq \mathcal{H}$ equipped with the supremum norm. Let $\mathcal{L}$ be the self-adjoint version of $(\partial_u + A^*)(BB^*)^{-1}(\partial_u - A)$ with boundary conditions $\omega(0) = 0$, $\omega'(1) = A\omega(1)$ on $\mathcal{H}$. From [14], Theorem 3.3, we know that the stationary distribution of the $\mathcal{H}$-valued SDE (3.1) coincides with $\nu$. By taking expectations on both sides of [14], equation (3.10) in the stationary state, we find that $m$ solves the boundary value problem

$$(\partial_u + A^*)(BB^*)^{-1}(\partial_u - A)m(u) = 0 \qquad \forall u \in (0,1),$$
$$m(0) = x^-, \qquad m'(1) = Am(1).$$

Define $U : E \to \mathbb{R}$ by $U(\omega) = \log(\psi(\omega))$ for all $\omega \in E$. We then have $d\mu = \exp(U(X))\,d\nu$ and the Fréchet derivative $F = U'$ is given by

$$\begin{aligned}(F(\omega - m))(u) = &-\nabla V(\omega(1))\delta(1 - u) \\ &+ D^2V(\omega(u))A\omega(u) + A^*\nabla V(\omega(u)) - \nabla \Phi(\omega(u))\end{aligned} \quad (6.3)$$

for all $\omega \in E + m$, where $\delta_1 \in E^*$ is a Dirac mass at $u = 1$ and $D^2V$ denotes the Hessian of $V$.

We check that the conditions of Theorem 3.4 are satisfied: from [14], Theorem 3.3, we know that $\mathcal{N}(0, -\mathcal{L}^{-1})$ is the distribution of $Z - m$. The only nontrivial point to be verified in assumption (A1) is the fact that $\mathcal{L}$ generates a contraction semigroup on $E$. This, however, follows immediately from the maximum principle. It follows from standard Sobolev embeddings that $\mathcal{H}^\alpha \subset E$ densely for every $\alpha > 1/4$ and [2], Lemma A.1 implies that $\mathcal{N}(0, -\mathcal{L}^{-2\alpha})$ is concentrated on $E$ in this case, so assumption (A2) also holds. Assumption (A3) is an immediate consequence of condition (M).

In order to check assumption (A4), define, for $n \geq 1$, the function $\delta_n(u) = n\chi_{[0,1/n]}$, where $\chi_A$ denotes the characteristic function of a set $A$. With this definition at hand, we define $F_n : E \to E$ by

$$\begin{aligned}(F_n(\omega - m))(u) = &-\nabla V(\omega(u))\delta_n(1 - u) \\ &+ D^2V(\omega(u))A\omega(u) + A^*\nabla V(\omega(u)) - \nabla \Phi(\omega(u)) \\ =: &(F_n^0(\omega))(u) + F^1(\omega(u)).\end{aligned}$$

Since $\mathcal{H}^\alpha$ is contained in some space of Hölder continuous functions (by Sobolev embedding), we have $\lim_{n\to\infty} \|F_n(\omega) - F(\omega)\|_{-\alpha} = 0$ for every $\omega \in E$. The locally uniform bounds on the $F_n$ as functions from $E$ to $\mathcal{H}^{-\alpha}$ follow immediately from condition (M), so it only remains to check the dissipativity bound (2.7).

We first use the representation (2.8) of the subdifferential in $E$ to check the condition $\langle \omega^*, F_n^0(\omega + y) \rangle \leq 0$ provided that $\|\omega\|_E$ is greater than some polynomially growing function of $\|y\|_E$. It follows from condition (M) and Hölder's inequality that there exists an increasing function $G:\mathbb{R}\to\mathbb{R}_+$ growing polynomially with $y$ such that

$$-\int_0^1 \langle \omega(u), \nabla V(\omega+y)\rangle \delta_n(1-u)|\omega^*|(du)$$
$$\leq -\int_0^1 (M(\omega(u),\ldots,\omega(u)) - G(|y(u)|))\delta_n(1-u)|\omega^*|(du)$$
$$\leq -\int_0^1 (c|\omega(u)|^{2p} - G(\|y\|_E))\delta_n(1-u)|\omega^*|(du)$$
$$= -(c\|\omega\|_E^{2p} - G(\|y\|))\int_0^1 \delta_n(1-u)|\omega^*|(du),$$

which is negative for $\|\omega\|_E$ sufficiently large. In order to check the corresponding condition for $F^1$, we treat the cases $p=1$ and $p>1$ separately, where $p$ is the exponent from condition (M).

In the case $p=1$, we can write $V(x)=\frac{1}{2}\langle x,Qx\rangle+\tilde V(x)$ for some positive definite matrix $Q$. We then have

$$F^1(x) = QAx + A^*Qx - QBB^*Qx + \bar F^1(x),$$

where $\bar F^1$ has sublinear growth at infinity. Condition (Q) then implies that there exists a constant $\gamma>0$ such that

$$\langle x, F^1(x)\rangle \leq -\gamma|x|^2 + |x||\bar F^1(x)|,$$

so that (2.7) follows from 5.4.

In the case $p>1$, it follows from condition (M) that

$$\langle x, F^1(x)\rangle = -\sum_i M(x,\ldots,x,Be_i)^2 + \bar F^1(x),$$

where $\bar F^1$ behaves like $o(|x|^{4p-2})$ at infinity. The nondegeneracy of $M$ thus implies that there exist constants $\gamma>0$ and $C$ such that

$$\langle x, F^1(x)\rangle = -\gamma|x|^2 + C,$$

so that (2.7) again follows from 5.4.

We finally check that $U$ is bounded from above. In the case $p>1$, this follows easily from condition (M). In any case, $V$ is bounded from below, so that in the case $p=1$, we have

$$U(\omega+m) \leq C + \int_0^1 (\langle M\omega(u), A\omega(u)\rangle - |B^*M\omega(u)|^2)\,du + \int_0^1 \tilde G(\omega(u))\,du$$

for some function $G$ behaving like $o(|x|^2)$ at infinity. It thus follows from condition (Q) that $U$ is indeed bounded from above. This concludes the verification of the assumptions of Theorem 3.4.

We now check that the assumptions of 4.5 hold. The only fact that remains to be checked is that (4.8) holds in our case. This can be verified easily since the nonlinearity is of the form

$$(F(\omega))(u) = G_1(\omega(u)) + G_2(\omega(1))\delta(1-u),$$

so it suffices to multiply the functions $G_i$ by smooth cut-off functions in order to get the required approximations of $F$.

From Theorem 2.6, we get that SDE (2.1) has a unique, mild solution for every initial condition $x_0 \in E$. Corollary 4.5 shows that the unique, ergodic invariant measure of SDE (2.1) is $\mu$. Converting from a Hilbert-space-valued SDE to an SPDE, as outlined in Remark 5.5, we find equation (8.3). This completes the proof of statement (a). Statement (b) follows directly from Theorem 4.10. □

REMARK 6.2. If $(BB^*)^{-1}A$ is symmetric, the matrix $A$ can be incorporated into the potential $V$ by choosing $A = 0$ and replacing $V(x)$ with $V(x) - \frac{1}{2}\langle x, (BB^*)^{-1}Ax\rangle$. In this case, the SPDE (6.1) simplifies to the more manageable expression

$$\partial_t x(t,u) = (BB^*)^{-1}\partial_u^2 x(t,u) - \nabla\Phi(x(t,u)) + \sqrt{2}\partial_t w(t,u)$$
$$\forall (t,u) \in (0,\infty) \times [0,1],$$
$$x(t,0) = x^-, \quad \partial_u x(t,1) = f(x(t,1)) \quad \forall t \in (0,\infty),$$
$$x(0,u) = x_0(u) \quad \forall u \in [0,1].$$

Similar simplifications are possible for the SPDEs considered in the remainder of this section and in the next.

Using the preconditioning technique described above, we can construct modified versions of the SPDE (6.1) which still have the same stationary distribution. In the preconditioned SDE (2.14), we take $\mathcal{G} = -\mathcal{L}^{-1}$, where $\mathcal{L}$ is the self-adjoint version of $(\partial_u + A^*)(BB^*)^{-1}(\partial_u - A)$ with boundary conditions $\omega(0) = 0$, $\omega'(1) = A\omega(1)$ on $L^2$.

THEOREM 6.3. *Let $A \in \mathbb{R}^{d\times d}$ be a matrix, let $B \in \mathbb{R}^{d\times d}$ be invertible, $f = -BB^*\nabla V$, assume that $V$ satisfies conditions* (M) *and* (Q) *and let $x^- \in \mathbb{R}^d$. Denote by $L$ the differential operator $(\partial_u + A^*)(BB^*)^{-1}(\partial_u - A)$ and consider the $\mathbb{R}^d$-valued SPDE*

(6.4a) $\quad \partial_t x(t,u) = -x(t,u) + y(t,u) + \sqrt{2}\partial_t \tilde{w}(t,u), \quad x(0,u) = x_0(u),$

*where $\tilde{w}$ is a $\mathcal{G}$-Wiener process, $\Phi$ is given by (5.3) and $y(t,\cdot)$ is the solution of the elliptic problem*

(6.5)
$$-Ly(t,u) = \nabla\Phi(x(t,u)) + A^*(BB^*)^{-1}f(x(t,u))$$
$$+ (Df)^*(x(t,u))(BB^*)^{-1}Ax(t,u)$$
$$y(t,0) = x^-, \quad \partial_u y(t,1) = Ay(t,1) + f(x(t,1)).$$

(a) *This SPDE has a unique, strong solution for every $x_0 \in \mathcal{C}_-([0,1], \mathbb{R}^d)$ and its stationary distribution coincides with the distribution of the solution of SDE (5.1).*

(b) *For every bounded function $\varphi : \mathcal{C}_-([0,1], \mathbb{R}^d) \to \mathbb{R}$ with bounded Fréchet derivative and every $x_0 \in \mathcal{C}_-([0,1], \mathbb{R}^d)$, we have*

$$\lim_{T \to \infty} \frac{1}{T} \int_0^T \varphi(x(t, \cdot)) \, dt = \mathbb{E}(\varphi(X)) \qquad \text{in probability,}$$

*where $X$ is the solution of (5.1).*

PROOF. Choose $\mathcal{H}$, $E$, $\mathcal{L}$, $m$, $\mu$ and $U$ as in the proof of Theorem 6.1. From [14], Theorem 3.3, we know that $\mathcal{G}$ is the covariance operator of the law of the solution of (5.2) and thus is positive definite, self-adjoint and trace class. We have already checked that (A1)–(A3) hold in the proof of Theorem 6.1. Furthermore, (A6)–(A7) are trivially satisfied for our choice of $\mathcal{G}$, so it only remains to check (A5) in order to apply Theorem 2.10. Note that the nonlinearity $F$ is of the form

$$(F(x))(u) = F_1(x(u)) + F_2(x(1))\delta(1-u)$$

for some functions $F_i : \mathbb{R}^n \to \mathbb{R}^n$. It follows from condition (M) that there exist constants $C$ and $N$ such that both of these functions satisfy

$$\langle x, F_i(x+y) \rangle \le C(1+|y|)^N$$

for every $x$ and $y$ in $\mathbb{R}^n$. The validity of (A5) follows immediately.

Applying Theorem 2.10, we obtain that SDE (2.14) has a unique, strong solution for every initial condition $x_0 \in E$. Corollary 4.9 shows that the unique, ergodic invariant measure of SDE (2.1) is $\mu$. Converting from a Hilbert-space-valued SDE to an SPDE, as outlined in Remark 5.6, we find equation (6.1). This completes the proof of statement (a). Statement (b) follows directly from Theorem 4.11. □

**7. Bridge path sampling.** In this section, we construct SPDEs which sample, in their stationary state, bridges from the SDE (5.1). That is, the stationary distribution of the SPDE coincides with the distribution of solutions of the SDE (5.1), conditioned on $X(1) = x^+$. The main results appear in Theorems 7.1 and 7.2.

Note that, for consistency with the other results in this paper, we construct an $E$-valued SPDE theory. However, this functional framework is not actually needed for this problem because the boundary conditions are linear; it is indeed possible in this case to use a Hilbert space theory. In that functional setting, results analogous to those in this section are mostly contained in [25] and [4]. We also refer to the monographs [5, 6] for related results.

Finally, note that the SPDE (7.1) was also derived (in the one-dimensional case) in [20].

We write $\mathcal{C}_-^+([0,1], \mathbb{R}^d)$ for the set of all continuous functions from $[0,1]$ to $\mathbb{R}^d$ which run from $x^-$ to $x^+$.

THEOREM 7.1. *Let $A \in \mathbb{R}^{d \times d}$ be a matrix, let $B \in \mathbb{R}^{d \times d}$ be invertible, $f = -BB^* \nabla V$, assume that $V$ satisfies conditions* (M) *and* (Q) *and let $x^-, x^+ \in \mathbb{R}^d$. Consider the $\mathbb{R}^d$-valued SPDE*

$$\partial_t x = (\partial_u + A^*)(BB^*)^{-1}(\partial_u - A)x - \nabla \Phi(x)$$
(7.1a)
$$\qquad - (Df)^*(x)(BB^*)^{-1} Ax - A^*(BB^*)^{-1} f(x) + \sqrt{2} \partial_t w,$$

(7.1b) $\quad x(t, 0) = x^-, \qquad x(t, 1) = x^+,$

(7.1c) $\quad x(0, u) = x_0(u),$

*where $\partial_t w$ is space-time white noise and $\Phi$ is given by (5.3).*

(a) *This SPDE has a unique, mild solution for every $x_0 \in \mathcal{C}_-^+([0,1], \mathbb{R}^d)$ and its stationary distribution coincides with the distribution of the solution of SDE (5.1), conditioned on $X(1) = x^+$.*

(b) *For every bounded, measurable function $\varphi : \mathcal{C}_-^+([0,1], \mathbb{R}^d) \to \mathbb{R}$ and every $x_0 \in \mathcal{C}_-^+([0,1], \mathbb{R}^d)$, we have*

$$\lim_{T \to \infty} \frac{1}{T} \int_0^T \varphi(x(t, \cdot)) \, dt = \mathbb{E}(\varphi(X) | X(1) = x^+) \qquad almost\ surely,$$

*where $X$ is the solution of (5.1).*

PROOF. Let $X$ and $Z$ be the solutions of the SDEs (5.1) and (5.2), respectively. From Lemma 5.2, we know that the density of the distribution $X$ with respect to the distribution of $Z$ is given by (6.2). Let $\mathcal{L}(Z | Z(1) = x^+)$ denote the conditional distribution of $Z$ and let $m : [0, 1] \to \mathbb{R}^d$ be the mean of this distribution. Using Lemma 5.3 and substitution, the density of $\mu = \mathcal{L}(X | X(1) = x^+) - m$ w.r.t. the centered distribution $\nu = \mathcal{L}(Z | Z(1) = x^+) - m$ is then given by

$$\psi(\omega - m) = c \cdot \exp\left( \int_0^1 \langle \nabla V(\omega(u)), A\omega(u) \rangle \, du - \int_0^1 \Phi(\omega(u)) \, du \right)$$

for all $\omega \in \mathcal{C}_-^+([0,1], \mathbb{R}^d)$ and some normalization constant $c$.

Consider the Hilbert space $\mathcal{H} = L^2([0,1], \mathbb{R}^d)$ and the embedded Banach space $E = \{\omega \in \mathcal{C}([0,1], \mathbb{R}^d) | \omega(0) = \omega(1) = 0\}$ equipped with the supremum norm. Define the operator $\mathcal{L}$ on $\mathcal{H}$ to be the self-adjoint version of $(\partial_u + A^*)(BB^*)^{-1}(\partial_u - A)$ with boundary conditions $\omega(0) = \omega(1) = 0$. From [14], Theorem 3.6, we know that the stationary distribution of the $\mathcal{H}$-valued SDE (3.1) coincides with $\nu$. By taking expectations on both sides of [14], equation (3.11) in the stationary state, we find that $m$ solves the boundary value problem

$$(\partial_u + A^*)(BB^*)^{-1}(\partial_u - A)m(u) = 0 \quad \forall u \in (0,1),$$
$$m(0) = x^-, \qquad m(1) = x^+.$$

Define $U: E \to \mathbb{R}$ by $U(\omega) = \log(\psi(\omega))$ for all $\omega \in E$. We then have $d\mu = \exp(U(\omega))\, d\nu$ and the Fréchet derivative $F = U'$ is given by

$$(7.2) \quad F(\omega - m) = D^2 V(\omega(u)) A\omega(u) + A^* \nabla V(\omega(u)) - \nabla \Phi(\omega(u))$$

for all $\omega - m \in E$.

Since (7.2) is the same as (6.3) without the terms involving delta functions, we can check that (A1)–(A4) hold in exactly the same way as in the proof of Theorem 6.1. From Theorem 2.6, we obtain that SDE (2.1) has a unique, mild solution for every initial condition $x_0 \in E$. Corollary 4.5 shows that the unique, ergodic invariant measure of SDE (2.1) is $\mu$. Converting from a Hilbert-space-valued SDE to an SPDE, as outlined in Remark 5.5, we find equation (8.3). This completes the proof of statement (a). Statement (b) follows directly from Theorem 4.10. □

Again, we study the corresponding result which is obtained from the preconditioned SDE (2.14). Since it is, in general, easier to invert the Laplacian with Dirichlet boundary conditions rather than $\mathcal{L}$, we choose $\mathcal{G} = -\mathcal{L}_0^{-1}$, where $\mathcal{L}_0$ is the self-adjoint version of $(BB^*)^{-1}\partial_u^2$ with boundary conditions $\omega(0) = \omega(1) = 0$ on $L^2$. This procedure leads to the following result.

THEOREM 7.2. *Let $A \in \mathbb{R}^{d \times d}$ be a matrix, let $B \in \mathbb{R}^{d \times d}$ be invertible, $f = -BB^* \nabla V$, assume that $V$ satisfies conditions* (M) *and* (Q) *and let $x^-, x^+ \in \mathbb{R}^d$. Consider the $\mathbb{R}^d$-valued SPDE*

$$(7.3a) \quad \partial_t x(t,u) = -x(t,u) + y(t,u) + \sqrt{2}\partial_t \tilde{w}(t,u), \qquad x(0,u) = x_0(u),$$

*where $\tilde{w}$ is a $\mathcal{G}$-Wiener process and $y(t, \cdot)$ is the solution of*

$$(7.4) \quad \begin{aligned} (BB^*)^{-1}\partial_u^2 y &= (BB^*)^{-1} A \partial_u x - A^*(BB^*)^{-1} \partial_u x \\ &\quad + A^*(BB^*)^{-1} A x + \nabla \Phi(x) \\ &\quad + (Df)^*(x)(BB^*)^{-1} A x + A^*(BB^*)^{-1} f(x), \end{aligned}$$
$$y(t,0) = x^-, \qquad y(t,1) = x^+,$$

*with $\Phi$ given by (5.3).*

(a) *This SPDE has a unique, strong solution for every $x_0 \in \mathcal{C}_-^+([0,1], \mathbb{R}^d)$ and its stationary distribution coincides with the conditional distribution of the solution of SDE (5.1), conditioned on $X(1) = x^+$.*

(b) *For every bounded function $\varphi: \mathcal{C}_-^+([0,1], \mathbb{R}^d) \to \mathbb{R}$ with bounded Fréchet derivative and every $x_0 \in \mathcal{C}_-^+([0,1], \mathbb{R}^d)$, we have*

$$\lim_{T \to \infty} \frac{1}{T} \int_0^T \varphi(x(t, \cdot))\,dt = \mathbb{E}(\varphi(X) | X(1) = x^+) \qquad \text{in probability,}$$

*where $X$ is the solution of (5.1).*

PROOF. The proof works in almost the same way as the proof of Theorem 6.3. The primary difference is that, in the present case, the operator $\mathcal{G}$ is not the inverse of $-\mathcal{L}$, but only of its leading order part. □

**8. Nonlinear filter/smoother.** Consider the $\mathbb{R}^d \times \mathbb{R}^m$-valued system of stochastic differential equations

(8.1)
$$dX = f(X)\,du + B_{11}\,dW^x, \qquad X(0) \sim \zeta,$$
$$dY = A_{21} X\,du + B_{22}\,dW^y, \qquad Y(0) = 0,$$

where $B_{11} \in \mathbb{R}^{d \times d}$, $A_{21} \in \mathbb{R}^{m \times d}$ and $B_{22} \in \mathbb{R}^{m \times m}$ are matrices, $(B_{11} B_{11}^*)^{-1} f$ is a gradient and $W^x$ (resp. $W^y$) is an independent standard Brownian motion in $\mathbb{R}^d$ (resp. $\mathbb{R}^m$). We will construct an SPDE which has the conditional distribution of $X$ given $Y$ as its stationary distribution. $\zeta$ is the density of the initial distribution for $X$ in (8.1). The main results are stated in Theorems 8.2 and 8.4.

REMARK 8.1. It is straightforward to extend the contents of this section to more general systems of the form

$$dX = (A_{11} X + B_{11} B_{11}^* \nabla V_1(X))\,du + B_{11}\,dW^x,$$
$$dY = (A_{21} X + A_{22} Y + B_{22} B_{22}^* \nabla V_2(Y))\,du + B_{22}\,dW^y$$

instead of equation (8.1). (We include the $A_{11} X$-term in the two previous sections.) We do not do so here because it would clutter the presentation.

THEOREM 8.2. *Let $A_{21} \in \mathbb{R}^{m \times d}$, $B_{11} \in \mathbb{R}^{d \times d}$ and $B_{22} \in \mathbb{R}^{m \times m}$ and assume that $B_{11}$ and $B_{22}$ are invertible. Let $f = -B_{11} B_{11}^* \nabla V$ and assume that $V$ satisfies conditions* (M) *and* (Q). *Let $\zeta$ be a $C^2$ probability density such that $\alpha = e^V \zeta$ satisfies*

(8.2) $$\max\{\log \alpha(x), \tfrac{1}{2} \langle \nabla \log \alpha(x), x \rangle\} \leq -\varepsilon |x|^2$$

*whenever $|x| \geq c$ for some constants $\varepsilon, c > 0$. Consider the $\mathbb{R}^d$-valued SPDE*

(8.3a)
$$\partial_t x(t, u) = (B_{11} B_{11}^*)^{-1} \partial_u^2 x(t, u) - \nabla \Phi(x(t, u))$$
$$+ A_{21}^* (B_{22} B_{22}^*)^{-1} \left( \frac{dY}{du}(u) - A_{21} x(t, u) \right) + \sqrt{2} \partial_t w(t, u)$$

*for all $(t,u) \in (0,\infty) \times [0,1]$ with boundary conditions*

(8.3b) $\quad \partial_u x(t,0) = -B_{11}B_{11}^* \nabla \log \alpha(x(t,0)), \qquad \partial_u x(t,1) = f(x(t,1))$

*for all $t \in (0,\infty)$ and initial condition*

(8.3c) $\qquad\qquad\qquad\qquad x(0,u) = x_0(u)$

*for all $u \in [0,1]$, where $\partial_t w$ is space-time white noise and $\Phi$ is given by (5.3).*

(a) *This SPDE has a unique, mild solution for every $x_0 \in \mathcal{C}([0,1], \mathbb{R}^d)$ and its stationary distribution coincides with the conditional distribution $\mu_{X|Y}$ of $X$ given $Y$ where $X,Y$ solve (8.1).*

(b) *For every bounded, measurable function $\varphi \colon \mathcal{C}([0,1], \mathbb{R}^d) \to \mathbb{R}$ and every $x_0 \in \mathcal{C}([0,1], \mathbb{R}^d)$, we have*

$$\lim_{T \to \infty} \frac{1}{T} \int_0^T \varphi(x(t,\cdot))\, dt = \mathbb{E}(\varphi(X)|Y) \qquad \text{almost surely,}$$

*where $X,Y$ solve (8.1).*

REMARK 8.3. The condition (8.2) on $\alpha$ seems to be quite artificial. On the other hand, if no a priori information is given on the distribution of $X(0)$, it is natural to assume that $X(0)$ is given by the invariant measure, in which case $\log \alpha = -V$, so that the assumptions on $\alpha$ are satisfied. In this case, the boundary conditions (8.3b) reduce to the more symmetric expression

$$\partial_u x(t,0) = -f(x(t,0)), \partial_u x(t,1) = f(x(t,1)).$$

PROOF OF THEOREM 8.2. Define $\bar{V} \colon \mathbb{R}^d \times \mathbb{R}^m \to \mathbb{R}$ by $\bar{V}(x,y) = V(x)$. We can then write the system $(X,Y)$ from (8.1) as an SDE of the form

(8.4) $\quad d\begin{pmatrix} X \\ Y \end{pmatrix} = A \begin{pmatrix} X \\ Y \end{pmatrix} du - BB^* \nabla \bar{V}(X,Y)\, du + B\, d\begin{pmatrix} W^x \\ W^y \end{pmatrix}$

with

$$A = \begin{pmatrix} 0 & 0 \\ A_{21} & 0 \end{pmatrix}, \qquad B = \begin{pmatrix} B_{11} & 0 \\ 0 & B_{22} \end{pmatrix}.$$

Let $(\bar{X}, \bar{Y})$ be the solution of the linear, $\mathbb{R}^d \times \mathbb{R}^m$-valued SDE

(8.5) $\quad d\begin{pmatrix} \bar{X} \\ \bar{Y} \end{pmatrix} = A \begin{pmatrix} \bar{X} \\ \bar{Y} \end{pmatrix} dt + B\, d\begin{pmatrix} W^x \\ W^y \end{pmatrix}$

with initial conditions

$$\bar{X}(0) \sim \mathcal{N}(0, \varepsilon^{-1}), \qquad \bar{Y}(0) = 0.$$

We can use Lemma 5.2 to obtain the density of the distribution $\mu_{XY}$ of $(X,Y)$ with respect to the distribution $\mu_{\bar{X}\bar{Y}}$ of $(\bar{X}, \bar{Y})$. Since the nonlinearity

$(f(X), 0)$ is orthogonal to the range of $B$ in $\mathbb{R}^d \times \mathbb{R}^m$, the resulting density is

$$\varphi(\omega, \eta) = \exp\left( V(\omega(0)) - V(\omega(1)) - \int_0^1 \Phi(\omega(u))\, du \right) \theta(\omega(0))$$

for all $(\omega, \eta) \in C([0,1], \mathbb{R}^d \times \mathbb{R}^m)$. Here, $\theta(\cdot)$ is the density of the distribution of $X(0)$, relative to the Gaussian measure $\mathcal{N}(0, \varepsilon^{-1})$. This density is proportional to $\exp(-V(x) + \log \alpha(x) + \frac{1}{2}\varepsilon|x|^2)$.

From [14], Lemma 4.4, we know that the conditional distribution $\mu_{\bar{X}|\bar{Y}}$ of $\bar{X}$ given $\bar{Y}$ exists and Lemma 5.3 shows that, since $\varphi$ does not depend on $Y$, $\mu_{X|Y}$ has density $\varphi$ with respect to $\mu_{\bar{X}|\bar{Y}}$. Let $m$ be the mean of $\mu_{\bar{X}|\bar{Y}}$. The density $\psi$ of $\mu = \mu_{X|Y} - m$ w.r.t. $\nu = \mu_{\bar{X}|\bar{Y}} - m$ is then given by

$$\psi(\omega - m) \propto \exp\left( \log \alpha(\omega(0)) + \frac{\varepsilon}{2}|\omega(0)|^2 - V(\omega(1)) - \int_0^1 \Phi(\omega(u))\, du \right)$$

for all $\omega \in C([0,1], \mathbb{R}^d)$.

Consider the Hilbert space $\mathcal{H} = L^2([0,1], \mathbb{R}^d)$ and the embedded Banach space $E = \mathcal{C}([0,1], \mathbb{R}^d) \subseteq \mathcal{H}$ equipped with the supremum norm. Define the formal second order differential operator

$$L = (B_{11}B_{11}^*)^{-1} \partial_u^2 x - A_{21}^* (B_{22}B_{22}^*)^{-1} A_{21} x.$$

Define the operator $\mathcal{L}$ to be the self-adjoint version of $L$ on $\mathcal{H}$ with boundary conditions $\omega'(0) = \varepsilon B_{11} B_{11}^* \omega(0)$ and $\omega'(1) = 0$. From [14], Theorem 4.1, we know that the stationary distribution of the $\mathcal{H}$-valued SDE (3.1) coincides with $\nu$. By taking expectations on both sides of [14], equation (4.2) in the stationary state, we find that $m$ solves the boundary value problem $-Lm(u) = A_{21}^*(B_{22}B_{22}^*)^{-1}\frac{dY}{du}(u)$ for all $u \in (0,1)$ with boundary conditions $m'(0) = \varepsilon B_{11}B_{11}^* m(0)$ and $m'(1) = 0$.

Define $U: E \to \mathbb{R}$ by $U(\omega) = \log(\psi(\omega))$ for all $\omega \in E$. We then have $d\mu = \exp(U(\omega))\, d\nu$. The Fréchet derivative $F = U'$ is given by

$$F(\omega - m) = \log \alpha(\omega(0)) \delta_0 + \varepsilon \omega(0) \delta_0 - \nabla V(\omega(1)) \delta_1 - \nabla \Phi(\omega(u))$$

for all $\omega \in E$, where $\delta_0, \delta_1 \in E^*$ are delta distributions located at 0, 1, respectively.

At this point, we are back in a situation that is very close to the one of Theorem 6.1 and we can check that (A1)–(A4) are satisfied. Note that (8.2) ensures that $U$ is bounded from above and that the term $\log \alpha(X(0))\delta_0$ appearing in $F$ satisfies (2.7). The various statements now follow from Theorem 2.6 and Corollary 4.5 as in Theorem 6.1. $\square$

In the preconditioned version of this theorem, we take $\mathcal{G} = -\mathcal{L}_0^{-1}$, where $\mathcal{L}_0$ is the self-adjoint extension of $(B_{11}B_{11}^*)^{-1}\partial_u^2$ with boundary conditions $\omega'(0) = 0$ and $\omega'(1) = \varepsilon B_{11} B_{11}^* \omega(1)$ for an $\varepsilon$ chosen so that (8.2) holds. This yields the following result, in which it is important to note that $\tilde{w}$ depends upon $\varepsilon$.

THEOREM 8.4. *Assume that the conditions of Theorem 8.2 hold and consider the $\mathbb{R}^d$-valued evolution equation*

$$\partial_t x(t, u) = -x(t, u) + y(t, u) + \sqrt{2} \partial_t \tilde{w}(t, u), \qquad x(0, u) = x_0(u),$$

*where $\tilde{w}$ is a $\mathcal{G}$-Wiener process and $y(t, \cdot)$ is the solution of the problem*

$$(B_{11} B_{11}^*)^{-1} \partial_u^2 y = A_{21}^* (B_{22} B_{22}^*)^{-1} \left( A_{21} x - \frac{dY}{du} \right) + \nabla \Phi(x),$$

*with boundary conditions*

$$\partial_u y(t, 0) = \varepsilon B_{11} B_{11}^* (y(t, 0) - x(t, 0)) - B_{11} B_{11}^* \nabla \log \alpha(x(t, 0)),$$
$$\partial_u y(t, 1) = f(x(t, 1)).$$

*As usual, $\Phi$ is given by (5.3).*

*(a) This SPDE has a unique, strong solution for every $x_0 \in \mathcal{C}([0,1], \mathbb{R}^d)$ and its stationary distribution coincides with the conditional distribution $\mu_{X|Y}$ of $X$ given $Y$, where $X, Y$ solve (8.1).*

*(b) For every bounded function $\varphi : \mathcal{C}([0,1], \mathbb{R}^d) \to \mathbb{R}$ with bounded Fréchet derivative and every $x_0 \in \mathcal{C}([0,1], \mathbb{R}^d)$, we have*

$$\lim_{T \to \infty} \frac{1}{T} \int_0^T \varphi(x(t, \cdot)) \, dt = \mathbb{E}(\varphi(X)|Y) \qquad \text{in probability,}$$

*where $X, Y$ solve (8.1).*

PROOF. The proof is very similar to that of Theorem 6.3, so we omit it. □

**9. Conclusions.** In this paper we derived a method to construct nonlinear SPDEs which have a prescribed measure as their stationary distribution. The fundamental relation between the drift of the SPDE and the density of the stationary distribution is in analogy to the finite dimensional case: if we augment the linear SPDE by adding an extra drift term of the form $F = U'$, the stationary distribution of the new SPDE has density $\exp(U)$ with respect to the stationary distribution of the linear equation.

Since the resulting SPDEs have unique invariant measures and are ergodic, they can be used as the basis for infinite dimensional MCMC methods. The applications in Sections 6, 7 and 8 illustrate this approach to sampling by constructing SPDEs which, in their stationary states, sample from the distributions of finite dimensional SDEs, conditioned on various types of observations. However, our analysis is limited to problems for which the drift is linear plus a gradient. Furthermore, in the case of signal processing, the analysis is limited to the case where the dependency of the observation on the signal is linear.

We have clear conjectures about how to generate the desired SPDEs in general, which we now outline. We start by considering the first two conditioning problems, 1 and 2. Since we consider the general nongradient case, the linear term can be absorbed into the nonlinearity and we consider the SDE

$$(9.1) \qquad dX = f(X)\,du + B\,dW^x, \qquad X(0) = x^-.$$

In the physics literature, it is common to think of the Gaussian measure induced by this equation when $f = 0$ as having density proportional to

$$q_{\mathrm{lin}}(Z) = \exp\left(-\int_0^1 \frac{1}{2}\left|\frac{dZ}{du}\right|_B^2 du\right).$$

If we denote by $\delta$ the variational derivative of path-space functionals such as this and consider the SPDE

$$\frac{\partial z}{\partial t} = \delta \ln q_{\mathrm{lin}}(z) + \sqrt{2}\frac{\partial W}{\partial t}$$

(the last term being space-time white noise), then this will sample from Wiener measure or Brownian bridge measure, depending on which boundary conditions are applied. This is an infinite-dimensional version of the Langevin equation commonly used in finite-dimensional sampling. General linear SPDEs derived similarly are proved to have the desired sampling properties in [14].

One can use the formal density $q$ given above, in combination with Lemma 5.2, to derive a formal density on path space for (9.1), proportional to

$$q_{\mathrm{non}}(X) = \exp\left(-\int_0^1 \frac{1}{2}\left|\frac{dX}{du} - f(X)\right|_B^2 + \frac{1}{2}\operatorname{div} f(X)\,du\right).$$

This density also appears in the physics literature and is known as the Onsager–Machlup functional [11]. The SPDEs, which we derived in Sections 6 and 7, may be found by considering SPDEs of the form

$$\frac{\partial x}{\partial t} = \delta \ln q_{\mathrm{non}}(x) + \sqrt{2}\frac{\partial W}{\partial t}$$

(the last term again being space-time white noise). Again, this may be seen as a Langevin equation. Calculating the variational derivative, we see that this SPDE has the form

$$(9.2) \qquad \frac{\partial x}{\partial t} = (BB^*)^{-1}\frac{\partial^2 x}{\partial u^2} - \Theta(x)\frac{\partial x}{\partial u} - \nabla_x \Phi(x) + \sqrt{2}\frac{\partial W}{\partial t},$$

where

$$\Theta(x) = (BB^*)^{-1} Df(x) - Df(x)^* (BB^*)^{-1}.$$

For bridge path sampling, the boundary conditions are those dictated by the bridging property. For free path sampling, the variational derivative includes a contribution from varying the right-hand end point, which is free, giving rise to a delta function; this leads to the nonlinear boundary condition.

When $f$ has a gradient structure, the operator $\Theta(x) \equiv 0$ and the SPDE is analyzed in this paper (in the case $A = 0$). When $\Theta(x) \neq 0$, it will, in general, be necessary to define a new solution concept for the SPDE in order to make sense of the product of $\Theta(x)$ with the derivative of $x$; in essence, we must define a spatial stochastic integral, when the heat semigroup is applied to this product term. The Stratonovich formulation of the density $q_{\text{non}}$ suggests that some form of Stratonovich integral is required. The case where the nongradient part of the vector field is linear is considered in this paper, and the provably correct SPDE in that case coincides with the conjectured SPDE above.

Turning now to the case of signal processing, we generalize the observation equation (1.2) to

$$dY = g(X, Y)\, dt + \tilde{B}\, dW^y, \qquad Y(0) = 0.$$

We can derive the appropriate SPDE for sampling in this case by utilizing the Onsager–Machlup functional above and applying Bayes' rule. Define

$$q_y(X, Y) = \exp\left(-\int_0^1 \frac{1}{2} \left|\frac{dY}{du} - g(X, Y)\right|^2_{\tilde{B}}\right).$$

The Onsager–Machlup density for sampling $(X, Y)$ jointly is

$$q(X, Y) := q_{\text{non}}(X) q_y(X, Y).$$

By Bayes' rule, the conditioned density for $X|Y$ will be proportional to $q(X, Y)$, with proportionality constant independent of $X$. Thus, the SPDE for sampling in this case should be

$$\frac{\partial x}{\partial t} = \delta \ln q(x, Y) + \sqrt{2} \frac{\partial W}{\partial t}$$

(the last term again being space-time white noise and $Y$ being the given observation). In the case where $g(X, Y)$ depends only on $X$ and is linear, and when $f(X)$ has the form considered in this paper, this SPDE is exactly that derived in this paper. Outside this regime, we are again confronted with the necessity of defining a new solution concept for the SPDE and, in particular, deriving a spatial stochastic integration theory. A related SPDE can also be derived when the observations are discrete in time. In this case, delta functions are introduced into the SPDE at the observation times; the theory introduced in this paper is able to handle this since a similar issue arises for the nonlinear boundary conditions considered here. Langevin SPDEs which solve the signal processing problem are discussed in [13]. In that paper, numerical experiments are presented which indicate that the conjectured SPDEs are indeed correct.

Finally, let us remark that we are currently unable to treat the case of multiplicative noise. We do not believe that this is a fundamental limitation of the method, but interpreting the resulting SPDEs will require much more careful analysis.

In addition to the extension of the Langevin equation to nongradient problems and more general observation equation, there are many other interesting mathematical questions remaining. These include the study of second order (in time $t$) Langevin equations, the development of infinite-dimensional hybrid Monte Carlo methods, the study of conditional sampling for hypoelliptic diffusions and the derivation of sampling processes for nonadditive noise.

M. HAIRER
A. M. STUART
J. VOSS
MATHEMATICS INSTITUTE
UNIVERSITY OF WARWICK
COVENTRY
WARWICKSHIRE CV4 7AL
UNITED KINGDOM
E-MAIL: hairer@maths.warwick.ac.uk